\documentclass{amsart}
\usepackage{ifthen}
\usepackage{subfigure}
\usepackage{amsmath}
\usepackage{amssymb}
\usepackage{euscript}
\usepackage[all]{xy}
\usepackage{varioref}
\usepackage{graphicx}
\usepackage{verbatim}
\usepackage{lscape}
\usepackage{xspace}
\usepackage{boxedminipage}
\usepackage{ifthen}
\newboolean{pdfOutput}
\setboolean{pdfOutput}{false}
\newboolean{abridged}
\setboolean{abridged}{false}
\ifthenelse{\boolean{pdfOutput}}
{
}
{
\usepackage{epsfig}
}
\newcommand{\drawfile}[1]
{
  \ifthenelse{\boolean{pdfOutput}}
  {
  \includegraphics{{#1}.pdf}%
  }
  {
  \epsfig{file={#1}.eps}%
  }
}
\theoremstyle{plain}
\newtheorem{theorem}{Theorem}
\newtheorem{corollary}{Corollary}
\newtheorem{lemma}{Lemma}
\newtheorem{proposition}{Proposition}

\theoremstyle{remark}
\newtheorem{definition}{Definition}
\newtheorem{example}{Example}
\newtheorem{remark}{Remark}

\newcounter{remarkCounter}
\newlength{\setBracketHeight}
\newcommand{\SetSuchThat}[2]{
  \settoheight{\setBracketHeight}{\ensuremath{#2}}
  \ensuremath{\left\{\left.{#1\rule{0cm}{\setBracketHeight}}\,
      \right|\,{#2}\right\}}}

\newcommand{\LieDer}{\ensuremath{\EuScript L}}
\newcommand{\hook}{\ensuremath{\mathbin{ \hbox{\vrule height1.4pt
        width4pt depth-1pt \vrule height4pt width0.4pt depth-1pt}}}}

\newcommand{\pd}[2]{\ensuremath{\frac{\partial{#1}}{\partial{#2}}}}
\newcommand{\R}[1]{\ensuremath{\mathbb{R}^{#1}}}
\newcommand{\C}[1]{\ensuremath{\mathbb{C}^{#1}}}

\newcommand{\Z}[1]{\ensuremath{\mathbb{Z}^{#1}}}
\newcommand{\Ha}[1]{\ensuremath{\mathbb{H}^{#1}}}
\newcommand{\Oc}[1]{\ensuremath{\mathbb{O}^{#1}}}
\newcommand{\CP}[1]{\ensuremath{\mathbb{CP}^{#1}}}
\newcommand{\RP}[1]{\ensuremath{\mathbb{RP}^{#1}}}
\newcommand{\HP}[1]{\ensuremath{\mathbb{HP}^{#1}}}
\newcommand{\OP}[1]{\ensuremath{\mathbb{OP}^{#1}}}
\newcommand{\Gr}[2]{\ensuremath{\operatorname{Gr}\left({#1},{#2}\right)}}
\newcommand{\Gro}[2]{\widetilde{\operatorname{Gr}}\left({#1},{#2}\right)}

\newcommand{\Spin}[1]{\ensuremath{\operatorname{Spin}\left({#1}\right)}}

\newcommand{\Cl}[1]{\ensuremath{\operatorname{C\ell}\left({#1}\right)}}
\newcommand{\GL}[1]{\operatorname{GL}\left({#1}\right)}

\newcommand{\slLie}[1]{\mathfrak{sl}\left({#1}\right)}

\newcommand{\Lm}[2]{\ensuremath{\Lambda^{#1} \left ( {#2} \right )}}

\newcommand{\Cohom}[2]{\ensuremath{ H^{#1} \left ( {#2} \right
    )}}
\newcommand{\Homol}[2]{\ensuremath{ H_{#1} \left ( {#2} \right
    )}}

 \DeclareMathOperator{\tr}{tr}
\renewcommand{\Re}{\operatorname{Re}}
\renewcommand{\Im}{\operatorname{Im}}

\DeclareMathOperator{\ind}{ind}
\newcommand{\OO}[1]{
  \ensuremath{
    \mathcal{O}
    \ifthenelse{\equal{#1}{0}}
      {}
      {\left({#1}\right)}
  }
}
\newcommand{\OOp}[2]{
  \ensuremath{
    \mathcal{O}
    \ifthenelse{\equal{#1}{0}}
      {}
      {\left({#1}\right)}
    \ifthenelse{\equal{#2}{1}}
      {}
      {^{\oplus{#2}}}
  }
}
\newcommand{\Bl}[2]{ 
  \ensuremath{
    \operatorname{Bl}_{#1}
    \left(
    {#2}
    \right)
  }
}

\newcommand{\extraSection}[2]
{
\ifthenelse{\boolean{abridged}}
  {
  }
  {
    \section{#1}
    \begin{center}
    \emph{This section will not be referred to
    subsequently, and may be skipped.}
    \end{center}
    \par{#2}
  }
}
\newcommand{\extraSubsection}[2]
{
\ifthenelse{\boolean{abridged}}
  {
  }
  {
    \subsection{#1}
    \begin{center}
    \emph{This subsection will not be referred to
    subsequently, and may be skipped.}
    \end{center}
    \par{#2}
  }
}
\newcommand{\extraStuff}[1]
{
\ifthenelse{\boolean{abridged}}
  {
  }
  {
    {#1}
  }
}
\newcommand{\Pts}{\ensuremath{P}}
\newcommand{\Lines}{\ensuremath{\Lambda}}
\newcommand{\Cor}{\ensuremath{F}}
\newcommand{\RT}[1]{\ensuremath{\check{#1}}}
\newcommand{\DRT}[1]{\ensuremath{\hat{#1}}}

\newcommand{\Epi}[2]{\ensuremath{\operatorname{Epi}\left({#1},{#2}\right)}}
\newcommand{\Algs}[1]{\ensuremath{\mathcal{A}_n}}
\newcommand{\divAlgs}[1]{\ensuremath{\mathcal{A}^+_n}}
\newcommand{\grLine}[3]{\ensuremath{\widetilde{\operatorname{Gr}}_{#1}
\left(#2,#3\right)}}
\def\cprime{$'$}
\newcommand{\Splus}{\ensuremath{S\Lambda^{2+}}}
\newcommand{\Sminus}{\ensuremath{S\Lambda^{2-}}}
\newcommand{\IN}[2]{\ensuremath{\#\left({#1} \cap {#2}\right)}}

\begin{document}
\title{Smooth projective planes}
\author{Benjamin McKay}
\address{University College Cork \\
  Cork, Ireland} \email{b.mckay@ucc.ie}
\date{\today} 
\thanks{Thanks to Robert Bryant and John Franks.}
\begin{abstract}
  Using symplectic topology and the Radon transform, we prove that
  smooth 4-dimensional projective planes are diffeomorphic to
  $\mathbb{CP}^2$. We define the notion of a plane curve in a smooth
  projective plane, show that plane curves in high dimensional regular
  planes are lines, prove that homeomorphisms preserving plane curves
  are smooth collineations, and prove a variety of results
  analogous to the theory of classical projective planes.
\end{abstract}
\thanks{This version of the paper incorporates minor corrections
  not included in the version published in \textbf{Geometriae Dedicata}.}
\maketitle
\tableofcontents
\section{Introduction}
A smooth projective plane is an object of topological
geometry, a manifold with a family
of submanifolds, called lines, satisfying the axioms
of projective plane geometry. A symplectic manifold
is an object of symplectic topology, a manifold $M^{2n}$
with a closed 2-form $\Omega$ for which $\Omega^n \ne 0$.
This article builds a bridge between symplectic topology and
topological geometry and
\begin{enumerate}
\item constructs differential forms on smooth projective planes via
  the Radon transform
\item proves that these forms tame lines (in the sense of Gromov's
  work on elliptic differential equations)
\item proves that smooth 4-dimensional projective planes are
  symplectomorphic to $\CP{2}$. (Previously only homeomorphism with
  $\CP{2}$ was known. Diffeomorphism of smooth projective planes
  of all other dimensions with $\RP{2},\HP{2}$ or $\OP{2}$
  was recently proven by Kramer \& Stolz.)
\item defines a concept of plane curve in a smooth projective plane
\item proves that the dual (i.e. set of tangent lines) of a plane
  curve is a plane curve
\item uncovers the lowest order differential invariant of a smooth
  projective plane (the tableau)
\item proves that the tableau determines the 2-jet of the smooth
  projective plane
\item proves that the tableau is a triality (in the sense of Cartan)
  just when the plane is regular (in the sense of Breitsprecher and
  B{\"o}di)
\item proves that regularity is invariant under projective duality
\item determines when a smooth projective plane agrees to second order
  with a classical projective plane
\item proves that on regular projective planes of dimension 4 or more,
  the differential equations defining plane curves are elliptic
  (determined for 4-dimensional projective planes, overdetermined for
  higher dimensions)
\item proves that the only plane curves in higher dimensional regular
  projective planes are lines
\item connects the analysis of plane curves in regular 4-dimensional
  projective planes to symplectic topology
\item proves that every regular 4-dimensional projective plane can be
  deformed through a family of 4-dimensional projective planes into
  one which is isomorphic to $\CP{2}$
\item proves that the differential system for plane curves in an
  irregular smooth 4-dimensional projective plane is a uniform limit
  of differential systems for plane curves in a regular projective
  plane.
\item proves that the differential system for plane curves in a
  regular projective plane determines the projective plane
\item characterizes the classical projective planes $\RP{2}$ and
  $\CP{2}$ in terms of local differential invariants
\item proves a collection of results about smooth quadrics in regular
  4-dimensional projective planes, showing that
\begin{enumerate}
\item they have a connected 10-dimensional moduli space
\item they are diffeomorphic to 2-spheres
\item through any 5 points, with no 3 colinear, there is a unique
  smooth quadric
\item the dual of a smooth quadric is a smooth quadric
\end{enumerate}
\item uncovers a dynamical system on the 2-torus which is the analogue
  of the elliptic curve found in the proof of Poncelet's porism.
\end{enumerate}
\section{Definition of smooth projective planes}
\begin{definition}
  An \emph{incidence geometry} is a triple
  $\left(\Pts,\Lines,\Cor\right)$ where $\Pts$ is a set (whose
  elements are called the \emph{points} of the geometry), $\Lines$ is
  a set (whose elements are called the \emph{lines} of the incidence
  geometry) and $\Cor \subset \Pts \times \Lines$ (whose elements are
  called the \emph{pointed lines} of the incidence geometry).  $\Cor$
  is called the \emph{correspondence space}, \emph{incidence
    correspondence} or \emph{flag space}.  We will say that a point $p
  \in \Pts$ is on a line $\lambda \in \Lambda$ if
  $\left(p,\lambda\right) \in \Cor$.  The \emph{dual} incidence
  geometry is $\left(\Pts^*,\Lines^*,\Cor^*\right)$ where
  $\Pts^*=\Lines$, $\Lines^*=\Pts$ and
  $\Cor^*=\SetSuchThat{\left(\lambda,p\right)}{\left(p,\lambda\right)\in\Cor}$.
  An incidence geometry is called a \emph{projective plane} if
\begin{enumerate}
\item any two distinct points $p_1,p_2 \in \Pts$ are on a unique line
  $p_1 p_2 \in \Lines$
\item any two distinct lines $\lambda_1,\lambda_2 \in \Lines$ have a
  unique common point $\lambda_1 \lambda_2 \in \Pts$ on them
\item there are at least 4 points, no three of which are on the same
  line.
\end{enumerate}
\end{definition}
\begin{definition}
  We call a projective plane \emph{topological} if $\Pts$ and $\Lines$
  are compact topological spaces, and the maps $p_1,p_2 \mapsto p_1
  p_2$ and $\lambda_1,\lambda_2 \mapsto \lambda_1 \lambda_2$ are
  continuous.
\end{definition}
\begin{definition}
  We call a topological projective plane \emph{smooth} if $\Pts$ and $\Lines$ are
  smooth manifolds, $\Cor
  \subset \Pts \times \Lines$ is a smooth embedded submanifold, and
  the maps $p_1,p_2 \mapsto p_1 p_2$ and $\lambda_1,\lambda_2 \mapsto
  \lambda_1 \lambda_2$ are smooth maps.
\end{definition}
Topological projective planes are studied by certain German
topologists; the standard reference is Salzmann et. al.
\cite{Salzmann:1995}.
\section{State of the art on smooth projective planes}
There is a useful characterization of smooth projective planes, due to
B{\"o}di and Immervoll:
\begin{definition}
  Suppose that we have a projective plane, for which $\Pts$ and
  $\Lines$ are closed smooth manifolds, both of dimensions $2n$ ($n \ge 0$ an
  integer), and suppose that $\Cor \subset \Pts \times \Lines$ is a
  $3n$-dimensional closed smoothly embedded submanifold, so that the
  canonical maps
\[
\xymatrix{
  & \Cor \ar[dr]^{\pi_{\Lines}} \ar[dl]_{\pi_{\Pts}} & \\
  \Pts & & \Lines }
\]
given by $\pi_{\Pts}(p,\lambda)=p$ and
$\pi_{\Lines}(p,\lambda)=\lambda$ are both submersions.  Call $\Pts$ a
\emph{smooth generalized plane}.  Consider the subsets
$\bar{\lambda}=\pi_{\Pts} \pi_{\Lines}^{-1}(\lambda) \subset \Pts$ and
$\bar{p}=\pi_{\Lines} \pi_{\Pts}^{-1}(p) \subset \Lines$ for $p \in
\Pts$ and $\lambda \in \Lambda$.  We will identify $\lambda \in
\Lambda$ with $\bar{\lambda}$, and still call it a line
($\bar{\lambda}$ is also often called a \emph{point row} in the
incidence geometry literature), while $\bar{p}$ will be called the
pencil through $p$. Lines are submanifolds of $\Pts$.
\end{definition}
\begin{theorem}[B{\"o}di \& Immervoll \cite{BodiImmervoll:2000}]\label{thm:BI}
  A smooth generalized plane is a smooth projective plane just when
  any two lines are transverse, and the pencils of any two points are
  transverse. Conversely, every smooth projective plane is a smooth
  generalized plane.
\end{theorem}
\begin{corollary}
Every submanifold of $\Pts \times \Lines$ which is $C^2$ close
enough to $\Cor$ defines a smooth projective plane; hence smooth
projective planes of dimension $2n$ depend on $n$ functions of $3n$
variables.
\end{corollary}
\begin{remark}
B\"odi \cite{Boedi:1997} remarks that it is not known if there are
real analytic smooth projective planes not isomorphic to the
standard real projective plane. It is easy to construct lots of real
analytic functions on products of projective planes, for example
eigenfunctions of the Laplacian. Thereby one can easily construct
lots of real analytic vector fields, and deform the flag space $F$
of the real projective plane. The deformations are too many to be
accounted for by the finite dimensional symmetry group. See Hilbert
\cite{Hilbert:1971} for proof that all symmetries of the standard
real projective plane are continuous, and that the group of
symmetries is the projective general linear group. Therefore there
are infinitely many real analytic smooth projective planes not
isomorphic to one another.
\end{remark}
\begin{theorem}[Freudenthal \cite{Freudenthal:1957}]
The dimension of a smooth projective plane is either 0, 2, 4, 8 or
16.
\end{theorem}
For proof, see Salzmann et. al. \cite{Salzmann:1995} p.258.  Ignoring
0, these happen to be the dimensions of the smooth projective planes
$\RP{2},\CP{2},\HP{2}$ and $\OP{2}$ respectively (see Salzmann et. al.
\cite{Salzmann:1995} for definitions); we call each of these four
spaces the \emph{model} of any smooth projective plane with the same
dimension. Zero dimensional projective planes are discrete, and
henceforth dimensions 0 and 2 will be largely ignored.
\begin{theorem}[Salzmann \cite{Salzmann:1967,Salzmann:1969},
\cite{Salzmann:1995} 51.29,
  L{\"o}wen \cite{Loewen:1995}, LeBrun \& Mason
  \cite{Lebrun/Mason:2002}] Two-dimensional smooth projective planes
  are diffeomorphic to the real projective plane.
\end{theorem}
A two-dimensional projective plane carries a projective
structure (a quite general type of path geometry, corresponding in
local coordinates to a second-order ordinary differential equation for
lines, see Bryant, Griffiths \& Hsu \cite{Bryant/Griffiths/Hsu:1995}).
Generic smooth projective plane structures on the real projective
plane are not projective connections (i.e. the lines are not geodesics
of any connection on the tangent bundle). LeBrun \& Mason were
concerned with projective connections, although their proof of this
theorem works for projective structures as well; most of their work in
that paper does not appear to apply to projective structures.
\begin{theorem}[Kramer \cite{Kramer:1994}]
  The space of points $\Pts$ of a smooth projective plane of positive
  dimension is homeomorphic to its model, as is the space of lines
  $\Lines$.
\end{theorem}
Kramer's theorem uses some quite deep differential topology. It is
much easier to show that each line is compact and connected, and that
the point space and line space are compact and connected; see Salzmann
et. al. \cite{Salzmann:1995} p. 225.  It is also not difficult to
show:
\begin{theorem}[Breitsprecher \cite{Salzmann:1995} pp. 257,262]\label{thm:Breit}
  The lines of a positive dimension smooth projective plane are
  homeomorphic to spheres, and the cohomology class of any line is the
  generator of the cohomology ring of $\Pts$.  The cohomology rings of
  $\Pts, \Lines$ and $\Cor$, and their relations under pullback, are
  identical to those of the model.
\end{theorem}
Henceforth, if the dimension is greater than 2, we will use this
identification of cohomology rings to orient $\Pts, \Lines, \Cor$ and
to orient each line and the pencil of each point, to match with the model.  In
particular, a pair of lines will intersect at a unique point, by
hypothesis, but moreover this intersection will be transverse by
theorem~\vref{thm:BI}, and in a projective plane of dimension 4 or
more the intersection will be positive by matching of cohomology.
\begin{definition}
  A \emph{morphism} of projective planes $\Pts_0 \to \Pts_1$ is a pair
  of maps $\Pts_0 \to \Pts_1$ and $\Lines_0 \to \Lines_1$, taking the
  flag space $\Cor_0 \subset \Pts_0 \times \Lines_0$ to the flag space
  $\Cor_1 \subset \Pts_1 \times \Lines_1$.  An isomorphism is often
  called a \emph{collineation}.
\end{definition}
\begin{theorem}[B{\"o}di \& Kramer \cite{Bodi/Kramer:1994}]\label{thm:BK}
  Every continuous isomorphism of smooth projective planes is smooth.
\end{theorem}
Note that the inverse is not assumed continuous, and it is not known
whether the inverse must be continuous.
\section{Affine charts}
Consider a smooth projective plane $\Pts$.  We will write down the
well-known affine charts on $\Pts$.  Pick any three points $0,X$ and
$Y \in \Pts$.  (Think of $0$ as the origin of an affine plane, and $X$
and $Y$ as the points at infinity on the $x$ and $y$ axes.)  Call
$\overline{0 X}, \overline{0 Y}$ and $\overline{X Y} \subset \Pts$ the
\emph{axes} associated to this choice of three points, and
$\overline{0},\overline{X}$ and $\overline{Y} \subset \Lines$ the
\emph{dual axes}. Define a map
\[
\alpha :p \in \Pts \backslash \overline{X Y} \mapsto \left( p_X,p_Y
\right) \in \left( \overline{0 X} \backslash X\right) \times \left(
  \overline{0 Y} \backslash Y\right)
\]
by
\begin{align*}
  p_X &= \left(p Y\right)\left(0 X\right) \\
  p_Y &= \left(p X\right)\left(0 Y\right).
\end{align*}
Define another map,
\[
\DRT{\alpha} : \lambda \in \Lines \backslash \overline{0} \mapsto
\left( \lambda_X,\lambda_Y \right) \in \left(\overline{0 X} \backslash
  0\right) \times \left(\overline{0 Y} \backslash 0\right) ,
\]
by
\begin{align*}
  \lambda_X &= \lambda \left(0 X\right), \\
  \lambda_Y &= \lambda \left(0 Y\right).
\end{align*}
\begin{lemma}
  $\alpha$ and $\DRT{\alpha}$ are diffeomorphisms.
\end{lemma}
\begin{proof}
  It is easy to check that if
\[
\alpha(p)=\left(p_X,p_Y\right)
\]
then
\[
p=\left(p_X Y\right)\left(p_Y X\right).
\]
Therefore $\alpha$ is smooth with smooth inverse, so a diffeomorphism.
Similarly, if
$\DRT{\alpha}(\lambda)=\left(\lambda_X,\lambda_Y\right)$, then \(
\lambda = \lambda_X \lambda_Y, \) so $\DRT{\alpha}$ is also a
diffeomorphism.
\end{proof}
$\DRT{\alpha}^{-1} \alpha$ is a diffeomorphism between $\Pts$ with
axes removed and $\Lines$ with dual axes removed.

A pointed line, i.e. an element $(p,\lambda) \in \Cor$, can be mapped
to $\overline{0 X} \times \overline{0 Y} \times \overline{0 X}$ in two
ways:
\begin{align*}
  \alpha(p,\lambda) &=
  \left(p_X,p_Y,\lambda_X \right) \\
  \DRT{\alpha}(p,\lambda) &= \left(p_X,\lambda_Y,\lambda_X \right).
\end{align*}
Given $\left(p_X,p_Y,\lambda_X\right)$, we can compute $p = \left(p_X
  Y\right) \left(p_Y X\right)$ (as before), and $\lambda =
p\lambda_X$.  It is easy to check that $\alpha$ and $\DRT{\alpha}$ are
diffeomorphisms
\begin{align*}
  \alpha &: \SetSuchThat{ \left(p,\lambda\right)\in \Cor} {p \notin
    \overline{OX} \cup \overline{XY} \text{ and } \lambda \ne OX} \to
  \left(\overline{0X} \backslash X\right) \times \left(\overline{0Y}
    \backslash \left\{0,Y\right\} \right) \times
  \overline{0X}. \\
  \DRT{\alpha} &: \SetSuchThat{(p,\lambda) \in \Cor}{\lambda \notin
    \bar{0} \cup \bar{Y}} \to \overline{0X} \times \left(\overline{0Y}
    \backslash \left\{0,Y\right\} \right) \times \left(\overline{0X}
    \backslash 0\right).
\end{align*}
Affine charts are defined on complements of submanifolds of the
obvious codimensions. Moreover, as in the model, these charts are
orientation preserving in projective planes of dimension 4 or more,
since at the intersection points of the various lines, the maps are
clearly orientation preserving, by positivity of intersection of
lines.
\begin{lemma}
  The lines and pencils of a smooth generalized plane are smooth
  embedded submanifolds. Lines meet transversely.
\end{lemma}
\begin{proof}
In the $\DRT{\alpha}$ chart on $\Cor$, this is immediate: we fix
$\lambda_X$ and $\lambda_Y$ and vary only $p_X$.  We can cover
$\Cor$ with three affine charts. An affine chart turns its axes into
obviously transverse submanifolds.
\end{proof}
\section{Blowup}
\begin{definition}
  Given a point $p_0 \in \Pts$, define the \emph{blowup}
  $\Bl{p_0}{\Pts}$ at $p_0$ to be the subset of $\Cor$ consisting of
  pairs $\left(p,\lambda\right) \in \Cor$ so that
  $\left(p_0,\lambda\right) \in \Cor$.
\end{definition}
\begin{lemma}
  The blowup of a smooth projective plane at a point is a smooth
  submanifold of the flag space.  Moreover the map
  $\left(p,\lambda\right) \in \Bl{p_0}{\Pts} \mapsto \lambda \in
  \bar{p}_0$ is a smooth fiber bundle. This fiber bundle admits the
  global section $\lambda \in \bar{p}_0 \mapsto
  \left(p_0,\lambda\right) \in \Bl{p_0}{\Pts}$; the image of this
  section is called the \emph{exceptional divisor} at $p_0$, and is
  also written $\bar{p}_0$. The map $\left(p,\lambda\right) \in
  \Bl{p_0}{\Pts} \mapsto p \in P$ is smooth, surjective, and a local
  diffeomorphism away from the exceptional divisor.
\end{lemma}
\begin{proof}
  In an affine chart $\DRT{\alpha}$, setting $X=p_0$, the open subset
  of the blowup intersecting the domain of that chart is given by
  points $(p,\lambda)$ with $p_Y = \lambda_Y \ne Y$ and $\lambda_X=X$
  and $p_X$ an arbitrary point of $\overline{OX}$ other than $X$.  We
  obtain a bijection
\[
\DRT{\alpha} : \SetSuchThat{ (p,\lambda) \in \Bl{p_0}{\Pts}} { \lambda
  \ne 0X,XY } \mapsto \left(p_X,p_Y\right) \in \overline{0X} \times
\left( \overline{0Y} \backslash \left\{0,Y\right\} \right ),
\]
(The map $\Bl{p_0}{\Pts} \to \bar{p}_0$ is $\left(p_X,p_Y\right)
\mapsto \lambda_Y=p_Y$, identifying an open set of the pencil
$\bar{p}_0$ with an open set of $\overline{OX}$.)  This covers all of
the blowup except for $\lambda \in \bar{O} \cup \bar{Y}.$

To cover those two pencils of lines, we switch coordinates, taking
$o=X,x=O,y=Y$.  Now we can try to use the other coordinates, $\alpha$
on $\Cor$. We then find that in terms of
$\left(p_x,p_y,\lambda_x\right)$, the blowup $\Bl{o}{\Pts}$ is given
by the equation $\lambda_x = o$. This gives a map
\[
\alpha : (p,\lambda) \mapsto \left(p_x,p_y\right)
\]
from
\[
\left\{ \lambda \ne ox,xy \text{ and } p \ne 0 \text{ and } p \ne
  \overline{xy} \right\} \subset \Bl{p_0}{\Pts}
\]
to
\[
\left( \overline{ox} \backslash \left\{o, x\right\} \right) \times
\left( \overline{oy} \backslash \left\{o,y \right\} \right).
\]
In this chart, the map $\Bl{p_0}{\Pts} \to \bar{p}_0$ is expressed as
$\left(p_x,p_y\right) \mapsto \lambda = po$ which we can map to the
$x$ or the $y$ axis.

Where both $\alpha$ and $\DRT{\alpha}$ are defined, we easily compute
$\left(p_x,p_y\right)$ in terms of $\left(p_X,p_Y\right)$ and vice
versa, via diffeomorphisms.  We have now covered all of the blowup
except for the points $\left(o,oy\right),\left(y,oy\right)$ and the
set of points of the form $\left(p,ox\right)$.

Swapping $x$ and $y$, so that the blowup is at $p_0=y$, and then
checking the smoothness of all of the maps, covers all of the blowup
except for
\[
\left(o,oy\right),\left(y,oy\right),
\left(o,ox\right),\left(x,ox\right).
\]
Finally, choosing any other choice of affine charts, perhaps
perturbing the $x$ and $y$ points slightly, covers the rest of the
blowup.
\end{proof}
\section{Hopf fibrations}
\begin{lemma}
Pick a point $p_0 \in \Pts$
and a line $\lambda_0 \in \Lines$
with $p_0$ not on $\lambda_0$. The map
\[
f : p \in \Pts \backslash p_0 \mapsto
\left(p p_0\right) \lambda_0 \in \bar{\lambda}_0
\]
is a fiber bundle mapping, with fiber
above $q \in \bar{\lambda}_0$
the punctured line $\overline{p_0 q} \backslash p_0$.
\end{lemma}
\begin{proof}
As in Salzmann et. al. \cite{Salzmann:1995} p. 252;
pick any $q \in \bar{\lambda}_0$, and
let $U = \bar{\lambda}_0 \backslash q,$
and let $A$ be the pencil of lines
through $q$, with $pq$ deleted.
Define the map
\[
\phi : \left(p,\lambda\right) \in U \times A
\mapsto
\left(p_0 p \right)\lambda
\in f^{-1} U.
\]
These maps trivialize our fiber bundle.
\end{proof}
\begin{lemma}
For each point $p_0 \in \Pts$,
the map
\[
f : p \in \Pts \backslash p_0 \mapsto p p_0 \in \bar{p}_0
\]
is a fiber bundle map with fiber through $\lambda$
being $\bar{\lambda} \backslash p_0$.
Call this map $f$ the \emph{Hopf fibration} at $p_0$.
\end{lemma}
\begin{proof}
This is essentially the same map.
\end{proof}
\section{The infinitesimal Hopf fibration}
We need to introduce an infinitesimal analogue of the Hopf fibration.
Given a point $p_0 \in \Pts$, start by constructing the pullback
vector bundle:
\[
\xymatrix{
  \tau \ar[d] \ar[r] & \ker \pi'_{\Lines} \ar[d] \\
  p_0 \times \bar{p}_0 \ar[d] \ar[r] & \Cor \ar[d] \\
  \bar{p}_0 \ar[r] & \Lines }
\]
so that the fiber of $\tau$ over a point $\lambda \in \bar{p}_0$ is
$T_{p_0} \bar{\lambda}$.  We can clearly map
\[
\tau \to T_{p_0} \Pts
\]
by inclusion.
\begin{lemma} The map $\tau \to T_{p_0} \Pts$
is a smooth bijection.
\end{lemma}
\begin{proof}
By transversality of lines, this map is an injection
away from the zero section.  Rescaling vectors by positive numbers
gives an injective smooth map
\[
S\tau \to ST_{p_0} \Pts
\]
from the bundle of spheres
\[
S\tau=\left(\tau \backslash 0\right)/\R{+},
\]
to the single sphere
\[
ST_{p_0} \Pts = \left(T_{p_0} \Pts \backslash 0 \right)/\R{+}.
\]
Each fiber of the bundle $S\tau$ is the sphere of the tangent plane of
a line, and is taken by the identity map to that same sphere of that
same tangent plane.  The map $\tau \to T_{p_0} \Pts$ is just the
derivative of the map $\Bl{p_0}{\Pts} \to \Pts$ restricted to the
submanifold $\bar{p}_0 \subset \Bl{p_0}{\Pts}$.
We need to see why $\tau \to T_{p_0} \Pts$ is onto.
Lets suppose that it misses some open set. Then
this open set, by rescaling, must contain an open cone.
Picture taking local coordinates with origin at the
point $p_0$. We can dilate these coordinates freely,
zooming in on the origin. As we do, the lines through
$p_0$ become flatter, with as many derivatives as we
like, and we can therefore approximate them uniformly
by their tangent planes at $p_0$. Therefore
picking any point $p$ lying in our cone (in
the given system of coordinates),
the line $pp_0$ enters that open cone. So $\tau \to T_{p_0} \Pts$
has dense image, and by rescaling $S \tau \to ST_{p_0} \Pts$
must as well. By compactness of $S \tau$,
the image must be closed. Therefore $S \tau \to ST_{p_0} \Pts$
is a smooth bijection, and so $\tau \to T_{p_0} \Pts$
is also a smooth bijection. By Sard's lemma, the
map $S\tau \to ST_{p_0} \Pts$ identifies the fundamental
classes in cohomology. Again, by Sard's lemma,
the inverse map $ST_{p_0} \Pts \to S\tau$ is
a diffeomorphism near a generic point.
\end{proof}
We call $ST_{p_0} \to S\tau \to \bar{p}_0$ the \emph{infinitesimal Hopf fibration}.
\begin{theorem}[B{\"o}di \cite{Boedi:1996}]\label{thm:BoediTwo}
  The map $S \tau \to ST_{p_0} \Pts$ is a smooth homeomorphism.
In particular, the infinitesimal Hopf fibration is a topological
sphere bundle, and a smooth submersion near a generic point.
\end{theorem}
\begin{proof}
  Smoothness of $S \tau \to ST_{p_0} \Pts$ is obvious,
  while homeomorphism is just the topological
  pigeonhole principle (see Salzmann et. al. \cite{Salzmann:1995} p.
  251); B{\"o}di's proof is different, employing the theory of
  microbundles.
\end{proof}
\begin{definition}
  We call the line $\lambda(\ell) \in \Lines$ tangent to a given real line
  $\ell \subset T_p \Pts$
  the \emph{magnification} of $\ell$.\footnote{The term
    \emph{magnification} is intended to remind the reader of
    \emph{complexification} of a real line to a complex line in
    $\CP{2}$.}
\end{definition}
\begin{lemma}
  Magnification is continuous, and smooth near a generic point.
\end{lemma}
\begin{proof}
  The magnification is the map taking a real line $\ell \subset
  T_{p_0} \Pts$ through the homeomorphism
  $T_{p_0} \Pts \to \tau$ and then through the vector bundle map $\tau
  \to \bar{p}_0$.
\end{proof}
\section{The tangential affine translation plane}
We have a new approach to defining the tangent plane of a smooth
projective plane.
\begin{definition}
  An \emph{affine plane} is a choice of sets $\Pts,\Lines,\Cor$ with
  map $\Cor \to \Pts \times \Lines$ so that, if we write $\pi_{\Pts} :
  \Cor \to \Pts$, $\pi_{\Lines} : \Cor \to \Lines$, as usual, and call
  the sets $\bar{\lambda}=\pi_{\Pts} \pi_{\Lines}^{-1}(\lambda)$
  \emph{lines} then
\begin{enumerate}
\item any two points $p_1,p_2$ lie on a unique line $p_1p_2$ and
\item for any point $p$ and line $\lambda$, there is a unique line
  $\lambda_p$ (called the \emph{parallel} to $\lambda$ through $p$) so
  that $p$ lies on $\lambda_p$ and either $\lambda_p=\lambda$ or
  $\lambda_p$ has no point in common with $\lambda$ and
\item there are three points not contained in any line.
\end{enumerate}
Following B{\"o}di \& Immervoll p. 66, we call an affine plane
\emph{smooth} (\emph{topological}) if
\begin{enumerate}
\item the maps $p_1,p_2 \mapsto p_1p_2$ and
  $p,\lambda\mapsto\lambda_p$ are smooth (continuous) and
\item the intersection point of two lines is unique, if it exists, and
  depends smoothly (continuously) on the choice of the lines, and
  exists for an open set of lines, and
\item there are four points with no three on a common line.
\end{enumerate}
A \emph{translation} of an affine plane is a map taking points to
points, lines to parallel lines, and preserving the flag space.  We
call an affine plane a (topological) [smooth] \emph{translation plane}
if the group of (homeomorphic) [diffeomorphic] translations acts
transitively on points.
\end{definition}
\begin{definition}
  Given a smooth projective plane $\Pts$ and a point $p_0 \in \Pts$,
  let $\Pts_0=T_{p_0} \Pts$, let $\tau \to \bar{p}_0$ be the vector
  bundle defined above, let $E_0=\bar{p}_0 \times P_0 \to \bar{p}_0$ be the
  trivial bundle, and let $\Lines_0$ be the quotient bundle $E_0/\tau
  \to \bar{p}_0$.  Consider the quotient map $Q : E_0 \to E_0/\tau$.
  Let $\Cor_0$ be the set of pairs $(p,l)$ in $E_0 \oplus_{\bar{p}_0}
  E_0/\tau$ for which $Q(p)=l$.  Define maps $(p,l) \in \Cor_0 \mapsto
  q \in \Lambda_0$ and $(p,l) \in \Cor_0 \mapsto p \in P_0$.
\end{definition}
\begin{lemma}\label{lemma:ATP}
  $\left(\Pts_0,\Lines_0,\Cor_0\right)$ is a topological affine
  translation plane with choice of origin $0 \in \Pts_0$, called the
  \emph{tangent plane} to $\Pts$ at $p_0$.  The translations are
  precisely the usual translations of $\Pts_0=T_{p_0} \Pts$, and are
  homeomorphisms.
\end{lemma}
\begin{proof}
  Clearly $Q$ is smooth and of constant rank, so that $\Cor_0$ is a
  smooth fiber subbundle of $E_0 \oplus _{\bar{p}_0} E_0/\tau$.  The
  maps $\pi_{\Pts}$ and $\pi_{\Lines}$ are obviously smooth
  submersions. The rest is proven in B{\"o}di \cite{Boedi:1997}.
\end{proof}
\begin{remark}
  For some smooth projective planes, the tangent plane is \emph{not}
  necessarily a
  smooth translation plane. The trouble comes from the smoothness of
  the map $p_1,p_2\mapsto p_1 p_2$. Indeed the map $p_2 \mapsto 0 p_2$
  is the magnification map.
\end{remark}
\section{The Radon transform}
\begin{definition}
  Let $\Pts$ be a smooth projective plane, with $\Lines$ its space of
  lines, and $\Cor$ its correspondence space. Take any top degree form
  $\eta$ on $\Lines$, and let $\RT{\eta}=\pi_{\Pts *} \pi_{\Lines}^*
  \eta$.  We call $\RT{\eta}$ the \emph{Radon transform} of $\eta$. If
  the dimension of $\Pts$ is $2n$, then $\RT{\eta}$ is an $n$-form.
\end{definition}
\begin{lemma}
  If $\eta$ is a volume form on $\Lines$ (i.e. a nowhere vanishing
  top degree form), then
  $\RT{\eta}$ is a closed form on $\Pts$, $\RT{\eta}^2$ is a volume
  form on $\Pts$, and $\phi^* \RT{\eta}$ pulls back to a positive
  volume form on any line.
\end{lemma}
\begin{proof}
  Suppose that $C\bar{\lambda}$ is a line through $p
  \in \Pts$, and $\RT{\eta}=0$ at $T_p \bar{\lambda}$.
  Using the map $\DRT{\alpha}$ on a dense open subset of
\[
\overline{0 X} \times \overline{0 Y} \backslash \left( 0 \times
  \overline{0 Y} \cup \overline{0 X} \times 0 \right),
\]
we write $\eta$ as
\[
\eta = f \, d \lambda_X \wedge d \lambda_Y,
\]
with $f>0$ and $d \lambda_X$ and $d \lambda_Y$ any
volume forms on the lines $\overline{0 X}$ and $\overline{0 Y}$
compatible with the orientations on those lines.
Pulling back,
\begin{align*}
  \pi_{\Lines}^* \eta
  &= f \, d \lambda_X \wedge d \lambda_Y \\
  &= f \, d \lambda_X \wedge \left( \pd{\lambda_Y}{p_X} \, dp_X +
    \pd{\lambda_Y}{p_Y} \, dp_Y \right).
\end{align*}
Pushing down,
\begin{align*}
  \RT{\eta} &=
  \pi_{\Pts *} \pi^*_{\Lines} \eta \\
  &= \left( \int f \pd{\lambda_Y}{p_X} \, d \lambda_X \right ) dp_X +
  \left( \int f \pd{\lambda_Y}{p_Y} \, d \lambda_X \right ) dp_Y.
\end{align*}
Therefore
\[
dp_X \wedge \RT{\eta} = \left( \int f \pd{\lambda_Y}{p_Y} \, d
  \lambda_X \right) \, dp_X \wedge dp_Y.
\]

We have to check signs: we need to ensure that we can consistently
orient lines to keep $\pd{\lambda_Y}{p_Y}>0$. The equations
\begin{align*}
  \lambda_Y &= \left( \left( p_X Y \right) \left( p_Y X \right )
    \lambda_X \right )
  \left( 0 Y \right) \\
  p_Y &= \left( 0 \lambda_Y \right) \left( \left(\lambda_X
      \lambda_Y\right) \left(p_X Y\right) \right),
\end{align*}
gives $\lambda_Y$ in terms of $p_Y$ and conversely, so these are
diffeomorphically mapped to one another, and so $\pd{\lambda_Y}{p_Y}
\ne 0$. By preservation of orientations under affine charts,
$\pd{\lambda_Y}{p_Y} > 0$.  As a consequence, $\RT{\eta} > 0$ on $T_0
\overline{0 X}$. We can pick the points $0$ and $X$ to be anywhere we
like, in particular pick $0=p$ and pick $X$ on the line
$\bar{\lambda}$. So $\RT{\eta}$ is a volume form on each line.

Given any point $p \in \Pts$, take two distinct lines $\lambda_1,
\lambda_2$ through $p$. They are transverse at $p$, and $\RT{\eta} \ne
0$ on each of their tangent planes at $p$, so $\RT{\eta}^2 \ne 0$ at
$p$.
\end{proof}
Note that a volume form $\eta$ on $\Lines$ exists
just when the dimension of $\Pts$ is 0,4,8 or 16 (i.e. not 2), as is
apparent from the cohomology.
\begin{corollary}
  If the dimension of a projective plane is 4, then the Radon
  transforms $\RT{\eta}$ of positive volume forms $\eta$ are all
  symplectomorphic, up to rescaling.
\end{corollary}
\begin{proof}
  Apply the Moser homotopy method.
\end{proof}
\begin{lemma}
  If $\Sigma \subset P$ is a compact oriented smooth submanifold
  (perhaps with boundary and corners) of dimension $n>0$ in a
  dimension $2n$ projective plane, then
\[
\int_{\Sigma} \RT{\eta} = \int_{\Lines} \IN{\Sigma}{\bar{\lambda}}
\, \eta
\]
where \( \IN{\Sigma}{\bar{\lambda}} \) is the number of
intersections of $\Sigma$ and $\bar{\lambda}$,
defined on the full measure subset
of $\lambda$ for
which $\Sigma$ and $\bar{\lambda}$ only intersect transversely.
\end{lemma}
\begin{remark}
We count \( \IN{\Sigma}{\bar{\lambda}} \)
keeping track of signs for positivity
or negativity of intersection. For a full measure
set of $\lambda \in \Lines$, the intersection
will be transverse, so we don't need to worry about
whether we are counting with multiplicity or not
at points of nontransverse intersection.
\end{remark}
\begin{proof}
  First, if we cut $\Sigma$ into two submanifolds $\Sigma_1$ and
  $\Sigma_2$, possibly with boundary and corners, which overlap only
  on their boundaries, then clearly
\[
\int_{\Sigma} \RT{\eta} = \int_{\Sigma_1} \RT{\eta} + \int_{\Sigma_2}
\RT{\eta}.
\]
On the other hand, the intersections of a line $\bar{\lambda}$ with
$\Sigma$ could occur either on $\Sigma_1$ or on $\Sigma_2$, but
possibly on both. However, the intersection $\Sigma' = \Sigma_1 \cap
\Sigma_2$ is a compact submanifold of dimension $n-1$, perhaps with
boundary and corners.  Define $\Cor'$ to be the pullback bundle
\[
\xymatrix{
  \Cor' \ar[d] \ar[r] & \Cor \ar[d] \\
  \Sigma' \ar[r] & \Pts.  }
\]
The lines striking $\Sigma'$ are the image of $\Cor' \to \Cor \to
\Lines$, so by Sard's theorem, counting dimension, the lines striking
$\Sigma'$ form a measure zero set. Therefore
\[
\int_{\Lines} \IN{\Sigma}{\bar{\lambda}} \, \eta =
\int_{\Lines} \IN{\Sigma_1}{\bar{\lambda}} \, \eta +
\int_{\Lines} \IN{\Sigma_2}{\bar{\lambda}} \, \eta.
\]
Therefore we only need to prove the result for ``small pieces'' of
$\Sigma$. The result holds for lines (for which it reduces to a
statement in cohomology), and therefore for any $\Sigma$ built out of
finitely many compact subsets of lines.

For any $\Sigma$, after cutting into enough submanifolds, we can see
that $\pi_{\Pts}^{-1} \Sigma$ is a smooth manifold with boundary and
corners, and $\pi_{\Pts}^{-1} \Sigma \to \Sigma$ is a fiber bundle,
with fiber over point $p$ the pencil $\bar{p}$.  Our integral is:
\[
\int_{\Sigma} \RT{\eta} = \int_{\pi_{\Pts}^{-1} \Sigma} \pi^*_{\Lines}
\eta.
\]
The map $\pi_{\Lines} : \pi_{\Pts}^{-1} \Sigma \to \Lines$ has
preimage at each point $\lambda$ given by all of the pairs
$(p,\lambda)$ with $p \in \Sigma$, i.e. $\Sigma \cap \bar{\lambda}$.
For $\Sigma$ a subset of a line $\bar{\lambda}_0$, these are
transverse positive intersections.  For generic $\Sigma$ and generic
$\lambda$ they are transverse intersections. The integrand vanishes
except at the points where $\pi_{\Lines}$ takes $\pi_{\Pts}^{-1}
\Sigma$ locally diffeomorphically to $\Lines$. Again by Sard's
theorem, the integral away from those points is just precisely the
integral of $\IN{\Sigma}{\bar{\lambda}} \, \eta$.
\end{proof}
\begin{corollary}
  The same is true if $\Sigma$ is a finite union of compact
  rectifiable submanifolds, possibly with boundaries and corners.
\end{corollary}
\begin{corollary}
  If $\Sigma \subset P$ is a compact $n$-dimensional rectifiable cycle
  (e.g. a line), then
\[
\int_{\Lines} \IN{\Sigma}{\bar{\lambda}} \, \eta = \int_{\Sigma}
\RT{\eta} = \left[\eta\right] \left[\Sigma\right]
\]
where the cohomology classes
$\left[\eta\right] \in \Cohom{2n}{\Lines,\R{}}=\R{}$ and
$\left[\Sigma\right] \in \Homol{2}{\Pts,\R{}}=\R{}$ are thought of as
numbers, using the standard basis for the cohomology.
\end{corollary}
\begin{corollary}
  The Radon transform $\eta \mapsto \RT{\eta}$ on projective planes of
  dimension 4 or more is injective.
\end{corollary}
\begin{theorem}\label{thm:CPTwo}
  Every smooth projective plane of dimension 4 is diffeomorphic to the
  complex projective plane.
\end{theorem}
\begin{proof}
  Such a projective plane is a compact 4-manifold with a symplectic
  structure.  Each line is diffeomorphic to a sphere, because we know
  its cohomology is that of a sphere. Moreover, its self-intersection
  is nonnegative, because in the orientation coming from the
  symplectic form, it belongs to a 4-dimensional family of spheres,
  and any two distinct spheres from that family intersect positively.
  Lalonde \& McDuff \cite{Lalonde/McDuff:1996} prove that a compact
  symplectic 4-manifold containing a symplectic sphere with
  nonnegative self-intersection is symplectomorphic to the complex
  projective plane, or a blowup of the complex projective plane at a
  finite number of points, or $S^2 \times S^2$. By cohomology, all of
  these are ruled out except the complex projective plane.
\end{proof}
\begin{remark}
  This idea is easy to generalize: Gromov \cite{Gromov:1985} p. 336
  suggests that a compact 4-manifold admitting a smooth incidence
  geometry of any reasonable type with compact embedded curves is
  diffeomorphic to $\CP{2}$, although he gives no details.  See
  B{\"o}di \& Immervoll \cite{BodiImmervoll:2000} for the definition
  of smooth incidence geometries.
\end{remark}
\begin{corollary}
  Every smooth projective plane of positive dimension is diffeomorphic
  to its model.
\end{corollary}
\begin{proof}
  We proved the result in dimension 4 above; for dimensions 8 and 16,
  the result was proven by Linus Kramer \& Stephan Stolz
  \cite{Kramer:2004}.
\end{proof}
\section{Immersed plane curves}
\begin{definition}
  A \emph{plane curve} in a smooth projective plane $\Pts$ of
  dimension $2n$ is a smooth immersion of manifolds $C^n \to \Pts$
  which is tangent to a line at each point.
\end{definition}
Note:
\begin{enumerate}
\item we do not ask the curve to be compact; e.g. in $\CP{2}$ this
  allows transcendental (i.e. nonalgebraic) curves
\item we do not allow singularities; picturing complex curves in
  $\CP{2}$, we would have to remove their singular points to fit this
  definition
\item in $\RP{2}$ (or any smooth projective plane of dimension 2) this
  definition allows all immersed curves, and is therefore useless. A
  reasonable concept of plane curve in a 2-dimensional projective
  plane, generalizing the concept of algebraic curve in $\RP{2}$, has
  never been formulated. There probably is one, in terms of Cartan's
  normal projective connection (see Bryant, Griffiths \& Hsu
  \cite{Bryant/Griffiths/Hsu:1995}).
\end{enumerate}
\begin{lemma}
Suppose that $P$ is a smooth projective plane, with space of lines $\Lines$.
The Radon transform of any volume form on $\Lines$ is
a positive volume form on every plane curve.
\end{lemma}
\begin{proof}
Plane curves are tangent to lines, and the Radon transform
of a volume form is positive on lines.
\end{proof}
\begin{lemma}
  In a smooth projective plane of dimension 4 or more, the homology
  class of any closed plane curve is a positive multiple
  of the homology class of a line.
\end{lemma}
\begin{proof}
  The homology is $\Homol{n}{\Pts}=\mathbb{Z}$, generated by
  $\left[\bar{\lambda}\right]$; see theorem~\vref{thm:Breit}.
Therefore $\left[C\right]=d \, \left[\bar{\lambda}\right]$,
for some integer $d.$ But
\[
0 < \int_C \RT{\eta} = d \int_{\bar{\lambda}} \RT{\eta}.
\]
\end{proof}
\begin{definition}
  The \emph{degree} of a plane curve $C$ in a smooth projective plane
  of dimension 4 or more is the ratio
\[
d=\frac{\left[C\right]}{\left[\bar{\lambda}\right]} \in \Z{+}.
\]
\end{definition}
\begin{lemma}\label{lemma:secant}
If $C$ is an embedded plane curve, and a sequence of points $p_j \in C$
approaches a limit $q \in C$, then the
secant lines $p_j q$ approach the tangent
line to $C$ at $q$.
\end{lemma}
\begin{proof}
There is a tangent line $\lambda$ to $C$ by definition,
and it is unique by transversality.
The points $p_j$ lie on the submanifold
$C$, and approach $q$. Therefore
in any local coordinates with
$q$ as origin, if $p_j$
is close enough to $q$, then
then $p_j$ lies near to $T_q C$.
Dilating coordinates as
needed, the line $p_j q$
is nearly a linear subspace in
some tiny coordinate ball.
Therefore $T_q \overline{p_j q} \to T_q C$
in the Grassmann bundle of linear subspaces.
The Gauss map $F \to \Gr{n}{TM}$ taking a
line to its tangent plane is clearly smooth,
and injective, and $F$ is a compact manifold. Therefore
the Gauss map is a topological embedding. So
the convergence of the tangent spaces
$T_{q}  \overline{p_j q} \to T_q C$ implies convergence
$p_j q \to \lambda$.
\end{proof}
\section{Polycontact system}
Define a field $\Theta$ of $2n$-planes on the
flag space $\Cor$ by assigning to each
point $\left(p,\lambda\right) \in \Cor$ the
$2n$-plane
\[
\Theta =
\ker \pi_P' \oplus
\ker \pi_{\Lambda}' \subset T \Cor.
\]
Call $\Theta$ the \emph{polycontact plane field}.
\begin{lemma}
$\Theta$ is a smooth plane field,
invariant under projective duality,
and
\[
\Theta_{\left(p,\lambda\right)}=
\left(\pi'_{\Pts}\right)
^{-1} T_p \bar{\lambda} \subset
T_{\left(p,\lambda\right)} \Cor.
\]
\end{lemma}
\begin{proof}
Smoothness is obvious, as is invariance
under projective duality. By definition,
$\bar{\lambda}=
\pi_{\Pts}\left(\pi^{-1}_{\Lines}\left(\lambda\right)\right)$.
Differentiating,
\begin{align*}
T_p \bar{\lambda}
&=
T_p \pi_{\Pts}\left(\pi^{-1}_{\Lines}\left(\lambda\right)\right) \\
&=
\pi_{\Pts}'(p,\lambda)
T_{(p,\lambda)} \pi^{-1}_{\Lines}\left(\lambda\right)\\
&=
\pi_{\Pts}'(p,\lambda)
\ker \pi'_{\Lines}(p,\lambda)\\
&=
\pi_{\Pts}'(p,\lambda)
\Theta_{(p,\lambda)},
\end{align*}
so
\[
\Theta_{(p,\lambda)} = \pi_{\Pts}'(p,\lambda)^{-1} T_p \bar{\lambda}.
\]
\end{proof}
\begin{definition}
Let $\pi : \Gro{n}{T\Pts} \to \Pts$ be the Grassmann bundle of
oriented $n$-planes in the tangent spaces of $\Pts$.
Let $\widetilde{\Theta} \subset T \Gro{n}{T \Pts}$
be the field of tautological planes,
\[
\widetilde{\Theta}_{\Pi} = \pi'(\Pi)^{-1} \Pi.
\]
\end{definition}
\begin{lemma}
The \emph{Gauss map}
\[
g : (p,\lambda) \in \Cor \to T_p \bar{\lambda} \in
\Gro{n}{T \Pts}
\]
is injective, and
\[
\Theta_{(p,\lambda)}=g'(p,\lambda)^{-1} \widetilde{\Theta}_{g(p,\lambda)}.
\]
In particular, if the Gauss map is an
immersion, then $\Theta$ is a subbundle
of $g^* \widetilde{\Theta}$.
\end{lemma}
The proof: unwind the definitions.
\section{Plane curves}
\begin{definition}
A \emph{generalized plane curve} in a projective
plane of dimension $2n$ is a Lipschitz map $\phi : C \to \Cor$
from an $n$-dimensional manifold $C$, perhaps
with boundary, so that
every differential form $\vartheta$ (of any degree)
on $\Cor$, vanishing on $\Theta$, pulls back to $\phi^* \vartheta = 0$.
A generalized plane curve is called \emph{basic} if
$\phi$ is injective on a dense
open set and intersects each fiber of $\pi_{\Pts} : \Cor \to \Pts$
on a discrete set of points.
\end{definition}
\begin{remark}
We will further generalize the notion of plane curve below
to allow singularities.
\end{remark}
\begin{lemma}
A continuously differentiable immersed
plane curve $\phi : C \to \Pts$
lifts to a continuous map $\Phi : C \to \Cor$ defined
by $\Phi(c)=\left(\phi(c),\lambda(c)\right)$, where
$T_{\phi(c)} \bar{\lambda}(c)=\phi'(c) T_c C$.
If the Gauss map is an
immersion, and $\phi$ is $C^{k+1}$, then
$\Phi$ is a $C^k$ basic generalized curve.
\end{lemma}
\begin{proof}
By hypothesis that $\phi$ is a plane curve,
$\Phi$ is defined; to see that $\Phi$ is
continuous, take any real immersed curve $c(t)$
on $C$, and look at its tangent lines
$\phi'(c(t)) c'(t)$, and magnify them:
$\lambda(c(t))=\lambda\left(\phi'(c(t)) c'(t)\right)$.
By theorem~\vref{thm:BoediTwo}, this map
is continuous.

Suppose that $\Cor \to \Gro{n}{T\Pts}$
is an immersion. Then
$\Phi$ maps to $\phi'$, so is
continuously differentiable.
To show that $\Phi : C \to \Cor$ is a
generalized plane curve, we will
show that $\Phi'(c) T_c C \subset \Theta$.
To see this, first note that $\pi_{\Pts} \Phi(c) = \phi(c)$.
Next, if $p=\phi(c)$ and $\lambda=\lambda(c)$, then
\begin{align*}
\Theta_{\left(p,\lambda\right)}
&=
\left(\pi'_{\Pts}\right)^{-1} T_p \bar{\lambda}
\\
&=
\left(\pi'_{\Pts}\right)^{-1} \phi'(c) T_c C
\\
&=
\left(\pi'_{\Pts}\right)^{-1} \left(\pi_{\Pts} \Phi\right)'(c) T_c C
\\
&=
\left(\pi'_{\Pts}\right)^{-1} \pi'_{\Pts} \Phi'(c) T_c C
\\
& \supset \Phi'(c) T_c C.
\end{align*}
\end{proof}
\begin{lemma}
A continuously differentiable
generalized plane curve $\Phi : C \to \Cor$ is the lift of
a plane curve $\phi : C \to \Pts$ just when
it is basic.
\end{lemma}
\begin{proof}
Define $\phi = \pi_{\Pts} \Phi$. Clearly
$\phi : C \to \Pts$ is an immersion because
$\Phi$ is not transverse to the fibers
of $\pi_{\Pts}$. We need to show that $\phi$ is
tangent to some line at every point. Write
$\Phi(c)=\left(\phi(c),\lambda(c)\right)$.
Since $\Phi'(c) T_c C \subset \Theta_{\Phi(c)}$,
\begin{align*}
\phi'(c) T_c C &=
\pi_{\Pts}' \Phi'(c) T_c C
\\
&\subset
\pi_{\Pts}' \Theta_{\left(\phi(c),\lambda(c)\right)}
\\
&=
\pi'_{\Pts} \left(\pi'_{\Pts}\right)^{-1} T_{\phi(c)} \bar{\lambda}(c) \\
&=
T_{\phi(c)} \bar{\lambda}(c).
\end{align*}
\end{proof}
\begin{definition}
If $\phi : C \to \Pts$ is an immersed projective
curve, then its \emph{dual curve} $\phi^* : C \to \Pts^*$
is $\phi^* = \pi_{\Lines} \Phi$, where $\Phi : C \to \Cor$
is the lift of $\phi$.
\end{definition}
\begin{corollary}\label{cor:dualCurves}
The dual curve $\phi^* : C \to \Pts^*$
of an immersed plane curve $\phi : C \to \Pts$
is an immersed plane curve in the dual plane, at every
point $c \in C$ where $\phi$ does not have
second-order tangency with a line.
\end{corollary}
\begin{lemma}
If $\eta$ is a positive volume form  on $\Lines$, and $\phi : C \to \Cor$
is a plane curve, then $\phi^* \pi^*_{P} \RT{\eta}$ is
a positive volume form on $C$.
\end{lemma}
\begin{lemma}
Generalized plane curves are precisely the
solutions of a smooth system of first order
partial differential equations.
\end{lemma}
\begin{proof}
The differential equations for
the lift $\Phi \subset \Cor$
are just $\vartheta=0$ for every
$\vartheta$ satisfying $\Theta=0$.
In local coordinates, a local
basis of such $\vartheta$ is easy
to write down, since $\Theta$ is
a vector bundle.
\end{proof}
\begin{lemma}
Plane curves (not generalized) are
precisely the solutions of a system of first order
partial differential equations.
\end{lemma}
\begin{proof}
Consider the map $(p,\lambda) \in \Cor \mapsto
T_p \bar{\lambda} \in \Gro{n}{T \Pts}$.
The requirement that an immersed manifold
be a plane curve is that its tangent
space lie in the image of this map, a
set of possible first derivatives in any
system of local coordinates.
\end{proof}
\section{Regularity}
Suppose that $\ell \subset T_{p_0} \Pts$ is a real line. Let
$\grLine{\ell}{n}{T_{p_0} \Pts}$
be the set of oriented $n$-planes in $T_{p_0} \Pts$
containing $\ell$. Our map $\Cor \to \Gro{n}{T \Pts}$ maps $\bar{p}_0
\to \Gro{n}{T_{p_0} \Pts}$.  Transversality of lines (which holds for
all smooth projective planes) is precisely the statement that
$\bar{p}_0$ can only intersect $\grLine{\ell}{n}{T_{p_0} \Pts}$
in at most one point.
\begin{definition}
  A smooth projective plane is \emph{regular} if for any point $p_0
  \in \Pts$ and any real line $\ell \subset T_{p_0} \Pts$, $\bar{p}_0$
  and $\grLine{\ell}{n}{T_{p_0} \Pts}$ inside $\Gro{n}{T_{p_0} \Pts}$ only have
  transverse intersections, if they have any intersections at all.
\end{definition}
\begin{example}
  The classical projective planes $\RP{2},\CP{2},\HP{2}$ and $\OP{2}$
  are all regular. Indeed by homogeneity of their collineation
  groups on their projectived tangent bundles, one need only
  check a single real line.
\end{example}
We will present several properties which are equivalent to regularity.
\begin{definition}\label{def:triality}
  A linear map $t : U \to V^* \otimes W$ (written $u \mapsto
  t_u$) with $\dim U = \dim V = \dim W$ is called a \emph{triality}
  (following Cartan \cite{Cartan:82}) if all of the linear maps
  $t_u$ are invertible, except at $u=0$.
\end{definition}
\begin{example}
  If $A$ is an algebra, then the multiplication map $A \otimes A \to
  A$ determines a map $t : A \to A^* \otimes A$, by $t_u=R_u$
  (right multiplication by $u \in A$). This map is a triality just
  when $A$ has no zero divisors. For $A=\R{},\C{},\Ha{}$ or $\Oc{}$,
  we will call this a \emph{classical triality}.
\end{example}
As is well known, the tangent planes to a Grassmannian
$\Gro{n}{\R{N}}$ at a point $E$ are intrinsically identified
with linear maps $E^* \otimes \left(\R{N}/E\right)$ as follows.
Consider the
  principal bundle
\[
\xymatrix{
  \GL{n,\R{}} \ar[r] & \Epi{\R{N}}{\R{n}} \ar[d] \\
  & \Gro{n}{\R{N}}.}
\]
The vertical map takes an epimorphism $T$ to its kernel.  The tangent plane to
$\Epi{\R{N}}{\R{n}}$ at any point is $\R{N*} \otimes \R{n}$.
Given a vector $\dot{E}$ tangent to the Grassmannian at $E$, take any
tangent vector $\dot{T}$ to $\Epi{\R{N}}{\R{n}}$ at $T$, so
that $\dot{T}$ maps to $\dot{E}$ under the derivative of the bundle
map.  Identifying $T_T \Epi{\R{N}}{\R{n}}$ with $\R{N*}
\otimes \R{n}$, consider $\left.T^{-1} \dot{T}\right|_{E} \in E^*
\otimes \left(\R{N}/E\right)$.  The reader can easily check
that $\dot{E} \mapsto \left.T^{-1} \dot{T}\right|_{E}$ smoothly
identifies the tangent vectors to the Grassmannian with linear maps.
The tangent plane to a subGrassmannian $\grLine{\ell}{n}{\R{N}}$
is precisely
the space of linear maps $\xi \in E^* \otimes \left(\R{N}/E\right)$
for which $\ell \subset \ker \xi$.
\begin{lemma}
  Map $\bar{p}_0$ to $\Gro{n}{T_{p_0} \Pts}$ by mapping a line
  $\lambda$ to its tangent plane $T_{p_0} \bar{\lambda}$. Regularity
  is just the requirement that this map is an immersion and that the
  induced maps on tangent planes of $\bar{p}_0$
\[
T_{\lambda} \bar{p}_0 \to E^* \otimes \left(T_{p_0} \Pts/E\right).
\]
are trialities.
\end{lemma}
\begin{proof}
Transversality to all of these subGrassmannians is precisely the
absence of kernel of all of the linear maps representing tangent
vectors.
\end{proof}
\begin{lemma}[Otte \cite{Otte:1992}]
  A smooth projective plane is regular just when the infinitesimal
  Hopf fibration is a smooth fiber bundle, i.e.  the tangent planes
  $T_{p_0} \bar{\lambda}$ of lines at a point $p_0$ divide the tangent
  plane $T_{p_0} \Pts$ into a fiber bundle (except at the origin).
  This is equivalent to smoothness of the magnification map $\ell
  \mapsto \lambda(\ell)$ (the line $\lambda$ so that $\ell \subset T_p
  \bar{\lambda}$, for $\ell \subset T_{p_0} \Pts$ a real line).
\end{lemma}
\begin{proof}
  Assume regularity. The map $\tau \to T_{p_0} \Pts$ above maps the
  tangent planes of lines injectively and smoothly into to the tangent
  plane. We need to show that it is a diffeomorphism away from the $0$
  section, and then we will use the fiber bundle structure of $\tau
  \backslash 0 \to \bar{p}_0$ to induce such a structure on $T_{p_0}
  \Pts \backslash 0$.  If the map $\tau \to T_{p_0} \Pts$ has a
  nonzero vector, say $v \in T_{p_0}$, in its image, then let $\ell$
  be the span of $v$.  Inside the Grassmannian, $\bar{p}_0$ strikes
  $\grLine{\ell}{n}{T_{p_0} \Pts}$
transversely at some point $\lambda$ (which is
  identified with $T_{p_0} \Pts \bar{\lambda}$ in $\Gro{n}{T_{p_0}
    \Pts}$).  Take $E \subset T_{p_0} \Pts$ any linear subspace for
  which $E \cap T_{p_0} \bar{\lambda} = \ell$.  For each real line
  $\ell'$ in $E$ close to $\ell$, the
  submanifold $\grLine{\ell'}{n}{T_{p_0} \Pts}$ must
  intersect $\bar{p}_0$ transversely at a single point, say
  $\lambda'$.  By compactness of $\bar{p}_0$, this must happen for all
  $\ell'$, not just those close to $\ell$.  Moreover, $\lambda'$
  clearly varies smoothly with $\ell'$ and $E$, by the implicit
  function theorem.  Therefore $\tau \to T_{p_0} \Pts$ has a smooth
  inverse.

  Conversely, assume that $\tau \to T_{p_0} \Pts$ is a diffeomorphism
  away from $0$, use the map $\tau \to \bar{p}_0$ to associate to each
  $\ell$ a smooth choice of $\lambda\left(\ell\right)$.  We need to show that
  for any nonzero vector $v \in T_{p_0} \bar{\lambda}$, and any vector
  $\xi \in T_{\lambda} \bar{p}_0$, $\xi(v) \ne 0$, thinking of $\xi$
  as a linear map. Suppose that $\xi(v)=0$ for some such $\xi \ne 0$
  and $v \ne 0$. Write $\xi$ as $\xi = T^{-1} \dot{T}$, so that
  $\dot{T} v = 0$ and $T v = 0$ and $v \ne 0$.  We must have $\xi =
  T^{-1} \dot T$ the velocity at $t=0$ of a curve $\lambda(t)$ in
  $\bar{p}_0$, and we can produce a curve $T(t)$ in $\Epi{T_{p_0}
    \Pts}{\R{n}}$ so that $\lambda(t)$ is its image in the
  Grassmannian. By assumption on $\xi$, $\dot{T}(0)v = T(0)v = 0$.
  Since $T_{p_0} \Pts$ fibers over $\bar{p}_0$ away from the origin,
  we can replace the single vector $v$ with a curve $v(t) \in T_{p_0}
  \Pts \backslash 0$ so that $\lambda(t)$ is its image in $\bar{p}_0$,
  $\dot{v}(0)$ is not in $T_{p_0} \bar{\lambda}$.  Therefore $v(t) \in
  T_{p_0} \bar{\lambda}(t)$ with $T_{p_0} \bar{\lambda}(t)$ the kernel
  of $T(t)$, i.e. $T(t) v(t)=0$. Putting our equations together,
\[
T(t) v(t) = \dot{T}(0) v(0) = 0.
\]
Plugging in Taylor expansions for $T(t)$ and $v(t)$ to these
equations, we find
\[
T(0) \dot{v}(0)=0,
\]
so that $\dot{v}(0) \in T_{p_0} \bar{\lambda}(0)$, contradicting our
hypothesis. Therefore the tangent planes of $\bar{p}_0$ are
trialities.
\end{proof}
\begin{corollary}
  Regularity occurs just when for any point $p_0 \in \Pts$ and real
  line $\ell \subset T_{p_0} \Pts$, $\bar{p}_0$ strikes
  $\grLine{\ell}{n}{T_{p_0} \Pts}$
  transversely at a single point, which depends
  smoothly on $p_0$ and $\ell$.
\end{corollary}
\begin{corollary}
  Regularity in the sense above is identical to regularity in the
  sense of Breitsprecher \cite{Breitsprecher:1971} and B{\"o}di
  \cite{Boedi:1996}.
\end{corollary}
\begin{proof}
  They define regularity as the property that the infinitesimal Hopf
  fibration is an isomorphism.
\end{proof}
\begin{corollary}\label{cor:generic}
  The generic point of the flag space of a smooth projective plane is
  regular.
\end{corollary}
\begin{proof}
  By B{\"o}di's theorem~\ref{thm:BoediTwo}, the infinitesimal Hopf
  fibration is a smooth homeomorphism. By Sard's lemma, the generic
  point of the target is a regular value. Therefore the generic point
  of the source is a regular point.
\end{proof}
\begin{theorem}
  Regularity is precisely smoothness of the tangent plane
  (i.e. the affine translation plane of lemma~\vref{lemma:ATP}).
\end{theorem}
\section{Trialities}
Recall from definition~\vref{def:triality} the concept of triality.
\begin{definition}
  A bilinear map $t : U \otimes V \to W$ is called a \emph{tableau}
  (see Bryant et. al. \cite{BCGGG:1991}).  Define the \emph{dual}
  tableau to be $t^* : V \otimes U \to W$, given by
  $t^*(v,u)=t(u,v)$.
\end{definition}
Consider a tableau $t : U \otimes V \to W$, with the property that
$t(u,v)\ne 0$ unless $u=0$ or $v=0$. This is identified with the
triality $t : U \to V^* \otimes W$ given by $t_u(v) = t(u,v)$.
\begin{remark}
  The concept of tableau does not require that $U,V,W$ have the same
  dimension. However, the concept of triality does.
\end{remark}
\begin{definition}
  Given a triality $t : U \to V^* \otimes W$, pick any elements
  $e_U \in U \backslash 0$ and $e_V \in V \backslash 0$.  Define
  $\epsilon_U : u \in U \mapsto t_u \left(e_V\right) \in V$ and
  $\epsilon_V : v \in V \mapsto t_{e_U}\left(v\right) \in W$.
  Given $v_0,v_1 \in V$ define
\[
v_0 v_1 = \epsilon_V^{-1} t\left(\epsilon_U^{-1}(v_0), v_1\right).
\]
This determines a real algebra with identity (not necessarily
commutative or associative).  We call this the \emph{triality algebra}
of these two points. We define the real part of any element $v$ of the
triality algebra to be
\[
\Re(v) = \frac{\tr\left(t_u\right)}{\dim U}.
\]
where $\epsilon_U(u)=v$.  We can think of this as a real number, or as
an element of the algebra. We define the imaginary part to be
$\Im(v)=v-\Re(v)$.
\end{definition}
\begin{remark}
  If $U,V,W$ have the same dimension, then the obvious notion of
  mapping tableau under linear transformations, taking a tableau
  $t$ and transformations $\left(g_U,g_V,g_W\right)$ with $g_U \in
  \GL{U}$, etc., and defining
\[
\left(g_U,g_V,g_W\right)t(u,v) = g_W t\left(g_U^{-1} u,
  g_V^{-1}v\right),
\]
is called \emph{isotopy of algebras} by algebraists (see Albert
\cite{Albert:1942}), perhaps an unfortunate term. The relevant objects
are really tableaux, not algebras.
\end{remark}
\begin{definition}
  Given an algebra $A$, and element $x \in A$, let $L_x : A \to A$ be
  the operation of left multiplication by $x$.  An algebra with
  identity element is called is called \emph{classical} if for any
  imaginary element $x \in A$, the endomorphism $L_x^2$ is a multiple
  of the identity element.
\end{definition}
\begin{lemma}
A triality algebra is classical just when the underlying triality is
classical (i.e. isomorphic to precisely one of the trialities of
right multiplication by real numbers, complex numbers, quaternions
or octave numbers).  Given a triality $U \subset V^* \otimes W$, the
triality algebra associated to a choice of two points $e_U \in U
\backslash 0$ and $e_V \in V \backslash 0$ is classical just when it
is classical at any pair of points.
\end{lemma}
\begin{proof}
  Write $L_x^2=-Q(x) \in \R{}$. Identifying $U$ and $V$ via
  $\epsilon_U$ and taking determinant
\[
\det\left(L_x^2\right)=Q(x)^n.
\]
Since $x \ne 0$, we must have $\det\left(L_x^2\right) =
\det\left(L_x\right)^2 \ne 0$ Therefore $Q$ is a definite quadratic
form.  If $Q$ is not positive definite, then pick some $x$ on which
$Q(x)=-1$, and compute
\[
\left(1-L_x\right)\left(1+L_x\right)=0.
\]
Because there are no zero divisors, $x=\pm 1$, is not imaginary. So
$Q$ is positive definite.  Therefore we can take $n$-th roots above
and get
\[
Q(x)=\det\left(L_x\right)^{2/n}.
\]
Let $V' \subset V$ be the set of imaginary elements. Then
\[
V' \hookrightarrow \slLie{V}
\]
extends to a representation of the Clifford algebra $\Cl{V',Q}$, since
it satisfies $L_x^2=-Q(x)$.  It is not a trivial representation, since
it is not trivial on $V'$.  Therefore it is a sum of irreducible
representations of $\Spin{n-1}$ in dimension $n$. Representation
theory tells us that such a representation exists only for $n=1,2,4$
or $8$, and is determined up to isomorphism.  In particular, it tells
us that the algebra is isomorphic to $\R{},\C{},\Ha{}$ or $\Oc{}$.
\end{proof}
\section{Regularity and the polycontact system}
\begin{definition}
  An \emph{adapted coframing} on the flag space $\Cor$ of a smooth
  projective plane is a choice of 1-forms
  $\vartheta^{\mu},\omega^{\mu},\pi^{\mu}$, $\mu=1,\dots,n$ on an open
  subset of $\Cor$ so that $\vartheta$ and $\omega$ are semibasic for
  $\Cor \to \Pts$, $\vartheta$ and $\pi$ are semibasic for $\Cor \to
  \Lines$ and $\Theta=\left(\vartheta=0\right)$.
\end{definition}
\begin{proposition}\label{prop:ZeroAdapted}
  Pick any adapted coframing, and calculate $d \vartheta = - \varpi
  \wedge \omega$ modulo $\vartheta$, where $\varpi$ is a combination
  of $\pi$ and $\omega$ 1-forms.  Then $\varpi$ has the form
  $\varpi^{\mu}_{\nu} = t^{\mu}_{\nu \sigma} \pi^{\sigma}$.  Call
  $t=\left(t^{\mu}_{\nu \sigma}\right)$ the \emph{tableau} of the
  adapted coframing.  The expression $t^{\mu}(u,v)=t^{\mu}_{\nu
    \sigma}u^{\nu} v^{\sigma}$ is a triality $t : \R{n} \otimes
  \R{n} \to \R{n}$ just where the projective plane is regular.
\end{proposition}
\begin{proof}
  We can ensure that $d \vartheta$ has the stated form, because
  unabsorbable $\omega \wedge \omega$ terms would prevent existence of
  integral manifolds, but we have the lines as integral manifolds, and
  similarly dually there can be no unabsorbable $\pi \wedge \pi$ terms
  (see Bryant et. al. \cite{BCGGG:1991}).

  Lets start off with a special choice of adapted coordinates. Take
  point $\left(p,\lambda\right) \in \Cor$, and take coordinates $x,y$
  on $\Pts$ near $p$, so that $p$ is $(x,y)=(0,0)$, and so that $T_p
  \bar{\lambda}$ is $dy=0$. Up on $\Cor$, we can let $\omega=dx$, and
  construct a suitable $\vartheta$ by demanding that $\vartheta = dy -
  z \, dx$, for some uniquely determined function $z$ defined on an
  open set.  Take $\pi^{\mu}$ arbitrary 1-forms semibasic for $\Cor
  \to \Lines$ and independent of $\vartheta$.

  I claim that for any vector $v \in T_{(p,\lambda)} \Cor$ tangent to
  the fiber $\bar{p}$, under motion with velocity $v$,
\[
\frac{dz^{\mu}_{\nu}}{dt} = t^{\mu}_{\nu \sigma} v^{\sigma},
\]
where we write $v^{\sigma}$ for $v \hook \pi^{\sigma}$.  To prove the
claim, extend $v$ to a vector field $X$ on an open subset of $\Cor$
and tangent to the fibers of $\Cor \to \Pts$, and compute $\LieDer_X
\vartheta = X \hook d\vartheta + d \left(X \hook \vartheta\right)$.
Note that since $X$ is tangent to the fibers, and $\vartheta$ is
semibasic, $X \hook \vartheta=0$.

Pick $\ell \subset (dy=0)$ any line, say spanned by a vector
$u=u^{\mu} \pd{}{x^{\mu}}$.  The open subset of the Grassmannian
$\grLine{\ell}{n}{T_{p_0} \Pts}$ on which $dx^1 \wedge \dots \wedge dx^n \ne 0$
consists of the planes $dy=z \, dx$ with $z^{\mu}_{\nu} u^{\nu} = 0$.
Therefore the tangent space of $\grLine{\ell}{n}{T_{p_0} \Pts}$
at $z=0$ is
$dz^{\mu}_{\nu} u^{\nu}=0$.  Consequently, the transversality claimed
is precisely expressed by the requirement that
\[
t^{\mu}_{\nu \sigma} v^{\sigma} u^{\nu} \ne 0
\]
for any $u$ and $v$.

Finally, if we change the choice of adapted coframing, say to
\[
\begin{pmatrix}
  \vartheta' \\
  \omega' \\
  \pi' \\
\end{pmatrix}
=
\begin{pmatrix}
  a & 0 & 0 \\
  b & c & 0 \\
  d & 0 & e
\end{pmatrix}
\begin{pmatrix}
  \vartheta \\
  \omega \\
  \pi
\end{pmatrix},
\]
then we can compute how $t$ changes:
\[
t'(u,v)=at\left(e^{-1}u,c^{-1}v\right).
\]
\end{proof}
\begin{corollary}
  The tableaux $t^{\mu}_{\nu \sigma}$ of the dual of a smooth
  projective plane are the dual tableaux of the original plane.
\end{corollary}
\begin{proof}
  A adapted coframing $\vartheta,\omega,\pi$ with $d \vartheta=-t
  \pi \wedge \omega$ for a smooth projective plane determines a
  adapted coframing
  $\vartheta^*=\vartheta,\omega^*=\pi,\pi^*=-\omega$ with $d
  \vartheta^* = -t^* \pi^* \wedge \omega^*$, $t^*$ the dual triality.
\end{proof}
\begin{corollary}
  The dual plane is regular just when the original plane is regular.
\end{corollary}
\section{Embedding the flag space into the
  Grassmannian bundles}
\begin{lemma}[Otte \cite{Otte:1992} 5.14]
  Write $\Gro{n}{T_p \Pts}$ for the Grassmannian of oriented
  $n$-planes in a tangent plane $T_p \Pts$, and $\Gro{n}{T \Pts}$ for
  the fiber bundle
\[
\xymatrix{
  \Gro{n}{T_p \Pts} \ar[r] & \Gro{n}{T \Pts} \ar[d] \\
  & \Pts.  }
\]
For a regular projective plane, the Gauss map
\[
\left(p,\lambda\right) \in \Cor \mapsto T_p \bar{\lambda} \in
\Gro{n}{T\Pts}
\]
is an embedding, and a fiber bundle mapping over $\Pts$.
\end{lemma}
\begin{remark}
  Projective duality gives an obvious dual to this lemma.
\end{remark}
\begin{remark}
  The reader might notice that $\Cor \to \Gro{n}{T \Pts}$ could
  perhaps be an embedding without regularity of the projective plane.
\end{remark}
\begin{proof}
  By transversality of lines, the map is injective.  To see that it is
  smooth, note that
\[
T_p \bar{\lambda} = \ker \pi'_{\Lines}\left(p,\lambda\right)
\]
is a vector bundle over $\Cor$. Transversality of the fibers $\bar{p}$
with the submanifolds $\grLine{\ell}{n}{T_p \Pts}$ forces the $\bar{p}$ to be
immersed submanifolds, because at each point $\lambda \in \bar{p}$, we
can pick any $\ell \subset T_{p} \bar{\lambda}$, and find $\bar{p}$
transverse to a submanifold of complementary dimension. Since the
fibers of $\Cor \to \Pts$ are compact (they are the $\bar{p}$
submanifolds), and the map $\Cor \to \Pts$ is a submersion, it is a
fiber bundle. Therefore since the fibers of $\Cor \to \Pts$ are
embedded into $\Gro{n}{T \Pts}$, it is clear from the local triviality
that $\Cor \to \Gro{n}{T \Pts}$ is an embedding.
\end{proof}
\section{Nondegenerate coordinates}
\begin{lemma}
  For any choice of point $p \in \Pts$ in a smooth projective plane,
  and choice of line $\lambda \in \Lines$, generic local coordinates
  $x,y : \Pts \to \R{n}$ on $\Pts$ defined near $p$, and generic local
  coordinates $X,Y : \Lines \to \R{n}$ defined near $\lambda$, the
  submanifold $\Cor \subset \Pts \times \Lines$ is given by equations
  in these coordinates, so that we can pick any three of $x,y,X,Y$,
  and the fourth will be a smooth function of those three.
\end{lemma}
\begin{proof}
  Pick any coordinates on $\Pts$ near $p_0$, and write them as $x,y$
  (with $x$ and $y$ each valued in $\R{n}$). Similarly write
  coordinates $X,Y$ on $\Lines$. After possibly a change of
  coordinates (which need only be a linear change of coordinates at
  worst), we can suppose that any choice of three of $x,y,X,Y$ gives
  coordinates on $\Cor$.  In particular, there must be functions
\begin{align*}
  x&=x(y,X,Y) \\
  y&=y(x,X,Y) \\
  X&=X(x,y,Y) \\
  Y&=Y(x,y,X).
\end{align*}
To check this, prove it first for linear functions $x,y,X,Y$ on vector
spaces $T_{\left(p_0,\lambda\right)} \Cor \subset T_{p_0} \Pts \times
T_{\lambda} \Lines$, and then the result is clear by the implicit
function theorem. Or just look at an affine chart.
\end{proof}
\begin{lemma}\label{lemma:regularityInCoordinates}
  A smooth projective plane is regular just when, in any nondegenerate
  coordinates, at any point $\left(x,y,X,Y\right) \in \Cor$, if we
  write $h^i_j$ for the inverse matrix of
\[
\pd{y^i}{Y^j},
\]
and write
\[
t^i_{jk} = \pd{^2 y^i}{x^j \partial X^k} - \pd{^2 y}{x^j \partial Y^l}
h^l_m \pd{y^m}{X^k}
\]
then for any nonzero vectors $\dot{x}$ and $\dot{X}$,
\[
t^i_{jk} \dot{x}^j \dot{X}^k \ne 0,
\]
i.e. $t^i_{jk}$ determines a triality.
\end{lemma}
\begin{proof}
  The tangent space at a point $(x,y,X,Y)$ to a line is given by the
  equation
\[
dy = \left.\pd{y}{x}\right|_{(x,X,Y)} \, dx.
\]
Parameterize the Grassmannian by associating to any plane $E \subset
T_{p_0} \Pts$ the matrix $M$ so that $dy = M \, dx$. Map $\Cor \to
\Gro{n}{T \Pts}$ by
\[
\left(x,X,Y\right) \mapsto \left(x,y,\pd{y}{x}\right).
\]
If $\ell \subset T_{p_0} \Pts$ is a real line contained in $T_{p_0}
\bar{\lambda}$, lets suppose that $\ell$ is spanned by a vector $a^i
\pd{}{x^i} + b^i \pd{}{y^i}$.  Without loss of generality, we can
suppose that $a^1=1$.  The subGrassmannian
$\grLine{\ell}{n}{T_{p_0} \Pts}$ is the set
of matrices $M$ so that $b = Ma$. To be transverse to this, we need
that whenever $b = \pd{y}{x} a$, for any curve $\lambda(t)$ in the
fiber $\bar{p}_0$, which passes through the given point with nonzero
velocity at $t=0$,
\[
\left.\frac{d}{dt}\right|_{t=0} \pd{y}{x} a \ne 0.
\]
Since we have to stay inside the pencil $\bar{p}_0$, we need
$\dot{x}=0$ and $\dot{y}=0$ on $\lambda(t)$.  Therefore
\[
0 = \dot{y} = \pd{y}{x} \dot{x} + \pd{y}{X} \dot{X} + \pd{y}{Y}
\dot{Y},
\]
from which we conclude that
\[
\dot{Y} = - \pd{Y}{y} \pd{y}{X} \dot{X}.
\]
Differentiating along such a curve gives precisely the stated
condition.
\end{proof}
\begin{lemma}
The flag space of any smooth projective plane is embedded in the
Grassmannian, i.e. the map $\Cor  \to \Gro{n}{T \Pts}$ is an
embedding, just when there is a vector $\dot{X}$ for which the
matrix $\dot{X}^k t^i_{jk}$ is not zero.
\end{lemma}
\begin{proof}
Immersion of the flag space in Grassmannian is expressed in
nondegenerate coordinates by requiring that $(x,X,Y) \mapsto
\pd{y}{x}$ be an immersion. However, it suffices for the fibers
$\bar{p}$ to be immersed in the Grassmannians, which is just the
requirement that every nonzero vector $\left(\dot{X},\dot{Y}\right)$
which is tangent to the fiber $\bar{p}$ must satisfy
\[
0 \ne \pd{^2 y}{x \partial X} \dot{X} + \pd{^2 y}{x \partial Y}
\dot{Y}.
\]
Since we must fix the point $(x,y)$, we find the embeddedness of the
flag space in the Grassmannian is precisely the condition that for
every vector $\dot{X}$, there is a vector $\dot{x}$ for which
$t^i_{jk} \dot{x}^j \dot{X}^k \ne 0$.
\end{proof}
We see clearly how regularity strengthens the requirement of being an
embedding. Note that $t^i_{jk}$ determines a triality: to each
$\dot{x}$ we associate the linear map $\dot{x}^k t^i_{kj}$.
\begin{corollary}
The flag space $\Cor$ of any smooth projective plane embeds into
$\Gro{n}{T \Pts}$ just when it embeds dually into $\Gro{n}{T
\Lines}$.
\end{corollary}
\section{Adapted coframings and nondegenerate coordinates}
\begin{lemma}
  Given any point of $\Cor$, there is a system of nondegenerate
  coordinates in which at any chosen point we can arrange
\[
\pd{Y}{x}=1, \pd{X}{y}=0, \pd{X}{Y}=0.
\]
The tableau $t$ is expressed at that point as
\[
t^{\mu}_{\nu \sigma} = \pd{^2 Y^{\mu}}{x^{\nu} \partial X^{\sigma}},
\]
in any adapted coframing which satisfies
$\vartheta=dy,\omega=dx,\pi=dX$ at that point.
\end{lemma}
\begin{proof}
  Start with any coordinates $z^I$ on $\Cor$ ($I=1,\dots,3n$), and
  nondegenerate coordinates $x,y$ on $\Pts$ and $X,Y$ on $\Lines$.
  Use the fact that the maps $z \mapsto (x,y)$ and $z \mapsto (X,Y)$
  have full rank, and the differentials of these maps have transverse
  kernels, to show that after perhaps a linear change of variables, we
  can arrange
\[
\pd{}{z}
\begin{pmatrix}
  x \\
  y \\
  X \\
  Y \\
\end{pmatrix}
=
\begin{pmatrix}
  0 & 1 & 0 \\
  1 & 0 & 0 \\
  0 & 0 & 1 \\
  1 & 0 & 0 \\
\end{pmatrix}.
\]
Now change variable to $z=\left(x,y,X\right)$, and you won't affect
these arrangements.  So we have $Y=Y(x,y,X)$ a function with the
required derivatives.

Write an adapted coframing $\vartheta,\omega,\pi$, with $\vartheta =
dy - p \, dx$ where
\[
p = - \left(\pd{Y}{y}\right)^{-1}\left(\pd{Y}{x}\right).
\]
We can arrange $\omega=dx$ and $\pi=dX$ at our distinguished point, by
changes of adapted coframing, so computing $d \vartheta$ gives
\[
t=\pd{Y}{x \partial X}
\]
as expected. Changing the coframing won't affect $t$, as long as it
doesn't change the coframe at that one point.
\end{proof}
\begin{lemma}
  The tableau $t=t^{\mu}_{\nu \sigma}$ of a smooth projective plane at
  a point of $\Cor$ determines, and is determined by the 2-jet of the
  map $\Cor \to \Pts \times \Lines$ at that point.
\end{lemma}
\begin{proof}
  The 1-jet of the map $\Cor \to \Pts \times \Lines$ determines and is
  determined by the $G$-structure whose sections are the adapted
  coframings. The 1-jet of that $G$-structure is precisely determined
  by the torsion (see, e.g., McKay \cite{McKay:unpub}), which in this
  case one easily computes to be the tableau.
\end{proof}
\section{Ellipticity of the differential equations for immersed
  plane curves}
\begin{lemma}
  On a regular projective plane of dimension 4 or more, the system of
  differential equations for plane curves is elliptic.
\end{lemma}
\begin{proof}
  Pick a basis of semibasic 1-forms $\vartheta^{\mu},\omega^{\mu}$ for
  the map $\Cor \to \Pts$ so that $\Theta = \left(\vartheta=0\right)$.
  Following Bryant et. al. \cite{BCGGG:1991}, ellipticity is the
  absence of a vector $v \in \Theta_{\left(p,\lambda\right)}$ so that
  $v \hook \varpi$ has rank 1.  As in the proof of
  proposition~\vref{prop:ZeroAdapted}, $v \hook \vartheta$ has nonzero
  determinant or vanishes, and therefore cannot have rank 1.
\end{proof}
Elliptic regularity results are quite involved for projective planes
of dimension $8$ and $16$, since the relevant equations are
overdetermined elliptic.  Nonetheless, with some effort one should be
able to identify in local coordinates a determined subsystem, and
prove elliptic regularity results for it. We are thereby encouraged
when studying plane curves to assume that they are smooth.

In principle, it seems possible that the equations
for immersed plane curves in some smooth irregular projective
plane might be elliptic, because $v \hook \vartheta$ might never
have rank 1, but might not still not have full rank. Examples of
irregular projective planes would be valuable.
\section{Cartan's count}\label{sec:count}
The reader who finds the material of this section mumbo-jumbo might
consult Bryant et. al. \cite{BCGGG:1991} for assistance.
\begin{lemma}
  The Cartan integers for the differential system of plane curves in a
  regular projective plane are $s_1=n,s_2=\dots=s_n=0$.
\end{lemma}
\begin{proof}
  The $v \hook \varpi$ matrices are all invertible, for $v \ne 0$, so
  they must all have a nonzero entry in the first column, and no
  linear combination of those entries vanish, so $\varpi$ has $n$
  independent 1-forms in its first column. There are only $n$
  independent directions available, modulo the independence condition
  (thinking of the differential system as a linear Pfaffian system),
  and therefore $s_2=\dots=s_n=0$.
\end{proof}
Therefore formally,\footnote{In stating that a system of smooth partial
differential equations has formal solution dependings on $s$ functions
of $d$ variables, we mean that the equations pass Cartan's test, with
last nonzero Cartan character $s_d=s$. In the real analytic category,
this will ensure that there is a well posed Cauchy problem for solutions
of the partial differential equations, with this generality. In the smooth
category, it tells us that while the corresponding Cauchy problem might
not really be well posed anymore, it is at least possible,
given $s$ functions of $d$ variables, to solve at a point in a formal Taylor
expansion solution which solves the equations at all orders.
Moreover, the formal solution will consist in Taylor coefficients
determined algebraically  by the Taylor coefficients of
those $s$ functions.
This is
often useful in trying to carry out approximation arguments which start
with a formal Taylor expansion solution, so even outside the real
analytic category, it is worthwhile to know the formal result of
Cartan's test. This is well but not widely known.}
plane curves depend on at most $n$ functions of
$1$ variable. However, we still have to check involutivity. Let us
compute the prolongation. Modulo $\omega=dx$, we can write
\[
\varpi^{\mu}_{\nu} = t^{\mu}_{\nu \sigma} \pi^{\sigma},
\]
with the $\pi^{\sigma}$ linearly independent.  Regularity will ensure
that $t^{\mu}_{\nu \sigma}$ is a triality. Plugging in $\pi = p
\omega$ to find the prolongation gives the equations
\[
t^{\mu}_{\nu \sigma} p^{\sigma}_{\tau} = t^{\mu}_{\tau \sigma}
p^{\sigma}_{\nu}.
\]
In terms of the algebra
\[
\left(xy\right)^{\mu} = t^{\mu}_{\nu \sigma} x^{\nu} y^{\sigma},
\]
this says that the prolongation consists in the $n \times n$ matrices
$p$ so that
\[
x(py)=y(px)
\]
for any $x$ and $y$.
\begin{definition}
  A linear transformation
  $p : A \to A$ of an algebra $A$ is
  an \emph{integral element} of the
  algebra if
\[
x(py) = y(px)
\]
for all $x, y \in A$.
More generally, an \emph{integral element} of a tableau $t : U
  \otimes V \to W$ is a linear map $p : U \to V$ so that
  $t\left(pu_0,u_1\right) =t\left(pu_1,u_0\right)$ for any $u_0,u_1
  \in U$ (see Bryant et. al. \cite{BCGGG:1991}).
\end{definition}
A simple calculation gives:
\begin{lemma}
  The integral elements of the algebra $\C{}$ are multiplications by
  complex numbers. More generally, any triality on $\R{2}$, after
  suitable linear transformation, has the form
\[
t^{\mu}_{\nu \sigma} = \delta^{\mu}_{\nu} \delta^0_{\sigma} +
a^{\mu}_{\nu} \delta^1_{\sigma}.
\]
The integral elements of the triality are the
matrices $p^{\mu}_{\nu}$ given by
\begin{align*}
  p^0_1 &= a^0_1 p^1_0 - a^0_0 p^1_1 \\
  p^0_0 &= a^1_1 p^1_0 - a^1_0 p^1_1.
\end{align*}
The set of such matrices is a 2-dimensional vector space.
\end{lemma}
\begin{corollary}
  The exterior differential system for plane curves in a 4-dimensional
  regular projective space is involutive; plane curves depend formally
  on 2 functions of 1 variable.  The system is elliptic. Formally
  (rigorously in the real analytic category), every connected real
  immersed curve sits in a unique maximal connected plane curve.
\end{corollary}
\begin{proof}
  See Bryant et. al. \cite{BCGGG:1991} for the relevant theory of
  noncharacteristic submanifolds for elliptic equations.
\end{proof}
Warning: a smooth but not analytic real immersed curve
does \emph{not} generally sit in any plane curve.
For instance, if we consider $\CP{2}$, take any real curve which
is contained in a complex curve. Now perturb the real curve
in some small open set, so that in that open set it no longer
lies on that complex curve, but in some other open set it still
does. Clearly there is no immersed complex curve containg the
real curve.
\begin{lemma}
  All integral elements vanish for all division algebras (i.e.
  trialities) on $\R{n}$ except in dimension $n=1,2$.
\end{lemma}
\begin{proof}
  By regularity, we can arrange with a simple change of coordinates at
  any required point that our algebra have a left identity (replace
  multiplication $uv$ with $u*v=L_{e_0}^{-1} (uv)$ where $L_{e_0}$ mean
  left multiplication by $e_0$; see Albert
  \cite{Albert:1942}).  It follows that $py=y(p1)$. So if we let
  $\epsilon=p1$, we find that
\[
x(y\epsilon)=y(x\epsilon)
\]
for all $x,y$. Let $R_{\epsilon}$ be the operation of right
multiplication by $\epsilon$. Define a new multiplication operation
$*$ by
\[
x*y = R_{\epsilon}^{-1}\left(x(y \epsilon)\right).
\]
The new multiplication is commutative and has no zero divisors, and
has the same identity element as the old one. By a famous theorem of
Heinz Hopf \cite{Hopf:1941} (see Springer \cite{Springer:1954} for a
beautiful proof, and also \cite{Hopf:1941,Behrend:1939}), the
dimension of a commutative (not necessarily associative) real algebra
without zero divisors must be 1 or 2.
\end{proof}
\begin{theorem}
  Let $\Pts$ be a regular projective plane of dimension 8 or 16.  The
  only basic plane curves in $\Pts$ are lines.
\end{theorem}
\begin{proof}
  Our differential system prolongs to a holonomic plane field, since
  the prolongation has dimension 0, and therefore the space $\Cor$ is
  foliated by the integral manifolds, a unique one through each point.
  But $\Cor$ is already foliated by (lifts of) lines, which
  are integral manifolds.
\end{proof}
\begin{theorem}
  For any smooth projective plane, $\left(\Pts,\Lines,\Cor\right)$, of
  dimension 8 or 16, and generic flag $(p,\lambda) \in \Cor$, there is
  no plane curve containing $p$ tangent to $\bar{\lambda}$, except
  $\bar{\lambda}$.
\end{theorem}
\begin{proof}
  See corollary~\vref{cor:generic}.
\end{proof}
\begin{remark}
  Note that this theorem is purely local: there are no ``little
  pieces'' of plane curves, not asking plane curves to be compact.  We
  can strengthen this slightly. Immersed plane curves are basic, but
  moreover every plane curve is basic on a dense open set unless it is
  a fiber of $\Cor \to \Pts$, i.e. a $\bar{p}$ pencil, which is a kind
  of degenerate plane curve, which we can just think of as a point. So
  roughly put, the only curves in high dimensional projective spaces
  are points and lines.
\end{remark}
\begin{remark}
  This theorem was previously unknown in any case, except for
\begin{enumerate}
\item the quaternionic projective plane (a folk theorem) and
\item the octave projective plane (Robert Bryant).
\end{enumerate}
Robert Bryant proved this result (but did not publish it) for the
octave projective plane using the canonical differential system on the
space of 8-planes calibrated by the $F_4$ invariant 8-form on the
octave projective plane. The same approach can be used (much more
easily) on the quaternionic projective plane. (It is also
easy to generalize to quaternionic projective spaces.) Generic smooth
projective planes of any dimension do not bear differential forms
calibrating their lines. Intuitively, the absence of calibrating forms
is explained by closed differential $n$-forms on a $2n$-manifold
depending on $\binom{2n}{n-1}$ functions of $2n$ variables, while a
smooth projective plane structure depends on $n$ functions of $3n$
variables. Small perturbations with compact support of a map $\Cor \to
\Pts \times \Lines$ which started out with all lines calibrated by a
given $n$-form will not remain calibrated by that $n$-form, and it is
easy to ensure that the lines will not be calibrated by anything.
\end{remark}
\begin{remark}
  It remains possible that a plane curve could exist on a smooth
  irregular projective plane of dimension 8 or 16, but the local
  nature of the proof would require the projective plane to have
  \emph{irregularities} at all points of the plane curve.  It should
  be possible to strengthen the results above to show that generic
  smooth projective planes of dimension 8 or 16 have no projective
  curves other than lines.
\end{remark}
\section{Foliating $\R{4}$ by 2-planes after Gluck \& Warner}
Following Gluck \& Warner \cite{Gluck/Warner:1983},
imagine a family of 2-planes in $\R{4}$ which foliate $\R{4}$ away
from 0. They prove that these 2-planes can be oriented continuously,
so that they all intersect positively,
and then determine a surface $\Sigma \subset \Gro{2}{\R{4}}$, in the
space of oriented 2-planes.  Pick a metric on $\R{4}$ and write each
2-plane $\Pi \subset \R{4}$ as $\xi \wedge \eta$ where $\xi$ and
$\eta$ are perpendicular unit-length 1-forms vanishing on $\Pi$.
Assume that $\Pi$ is oriented, so that the pair $\xi,\eta$ are
well-defined up to rotation.  A vector $v$ belongs to $\Pi$ just when
$v \hook \xi \wedge \eta=0$.  Given two 2-planes $\Pi,\Pi'$ we write
them as $\xi \wedge \eta,\xi'\wedge\eta'$, and clearly $0=\xi \wedge
\eta \wedge \xi' \wedge \eta'$ just when $\Pi \cap \Pi' \ne 0$.
Moreover, $\xi \wedge \eta \wedge \xi' \wedge \eta'>0$
just when the 2-planes have positive intersection.

Split $\xi \wedge \eta = \sigma_+ + \sigma_-$, where $\sigma_+$ is a
self-dual 2-form, and $\sigma_-$ is anti-self-dual.  If we write out
an orthonormal basis $\sigma_+^i$ for the self-dual 2-forms, and an
orthonormal basis $\sigma_-^i$ for the anti-self-dual 2-forms, then
$\sigma_{\pm}^i \wedge \sigma_{\pm}^j =\pm dV$ and $\sigma^i_+ \wedge
\sigma^j_- = 0$. Therefore writing
\[
\xi \wedge \eta = X_i \sigma_+^i + Y_i \sigma_-^i,
\]
we find
\[
1 = \frac{\left(\xi \wedge \eta\right)^2}{dV} = X_i^2 - Y_i^2,
\]
so that $\sigma_+$ and $\sigma_-$ belong to the unit spheres $S^+$ and
$S^-$ of self-dual and anti-self-dual 2-forms. Gluck \& Warner show
that map
\[
\left(\sigma_+,\sigma_-\right) : \Sigma \to S^+ \times S^-
\]
of a surface to the Grassmannian satisfies
$\left|\dot{\sigma}_-\right|> \left|\dot{\sigma}_+\right|$ just when
the family of 2-planes it represents locally smoothly foliates some
region in $\R{4}$ with 2-planes.  Consequently, if the entirety of
$\R{4}$ is foliated by 2-planes, then the image of this map $\Sigma
\to S^+ \times S^-$ is the graph of a smooth strictly contracting map
$S^- \to S^+$, and this strictly contracting map determines the
foliation.

All of this applies immediately to the tangent planes
to the lines through a point $p_0 \in \Pts$ in a 4-dimensional
smooth projective projective plane. These tangent planes
foliate $\R{4}=T_{p_0} P$ just when $p_0$ is a regular point.

A self-dual 2-form in the formalism is represented by a choice of
point $\omega_+ \in S^+$, and it is positive on a 2-plane $\xi \wedge
\eta$ just when $\xi \wedge \eta \wedge \omega_+ > 0$, a positive
volume form.
\section{Regular 4-dimensional projective planes}
McKay \cite{McKay:2003} develops the general theory of pseudocomplex
structures. In local coordinates, these are determined
systems of first order elliptic partial differential equations
for two functions of two variables. Globally, they are
a choice of such equations on each coordinate chart
of a 4-manifold, and having the same local solutions on overlaps
of coordinate charts.
\begin{theorem}\label{thm:taubes}
  The category of regular 4-dimensional projective planes is
  isomorphic to the category of compact, symplectically tameable
  pseudocomplex 4-manifolds which contain a pseudoholomorphic sphere
  with nonnegative selfintersection, which is isomorphic to the
  category of pseudocomplex structures on $\CP{2}$ tameable by a
  symplectic structure. In particular, smooth families from one
  category are smoothly equivalent to smooth families from the other.
  Every regular 4-dimensional smooth projective plane can be deformed
  through regular 4-dimensional projective planes tamed by a fixed
  symplectic form into a classical 4-dimensional projective plane.
\end{theorem}
\begin{remark}
  We will be brief in our analysis, for which all details are worked
  out completely in Lalonde \& McDuff \cite{Lalonde/McDuff:1996} and
  McKay \cite{McKay:2003}.
\end{remark}
\begin{proof}
  Given a regular 4-dimensional projective plane $\Pts$, pick a volume
  form $\eta$ on $\Lines$, in the cohomology class dual to
  $\left[\Lines\right]$, and let $\RT{\eta}$ be the Radon transform.
  Take the pseudocomplex structure to be $\Cor \subset \Gro{n}{T
    \Pts}$.  Pseudoholomorphic curves are precisely plane curves.  It
  follows by transversality of lines and intersection theory from
  \cite{McKay:2003} that $\Lambda$ is the moduli space of
  pseudoholomorphic spheres in the generating homology class.
  More specifically, if we have any compact pseudoholomorphic
  curve $\Sigma$ in the homology class of a line, and
  $p \in \Sigma$ is a smooth point of $\Sigma$, then
  take $\lambda$ the tangent line to $\Sigma$ at $p$.
  The surfaces $\bar{\lambda}$ and $\Sigma$ must have
  intersection number at least 2 at $p$ unless they
  are equal (see McKay \cite{McKay:2003}
  for proof). By cohomology calculations, $\Sigma=\lambda$.
  So the  map of categories is defined and injective.

  Conversely, suppose that $\Pts$ is a compact 4-manifold, has a
  symplectically tame pseudocomplex structure, and contains a
  pseudoholomorphic sphere with nonnegative selfintersection.  By
  results of McDuff (see Lalonde \& McDuff
  \cite{Lalonde/McDuff:1996}), this forces $\Pts$ to be
  symplectomorphic to $\CP{2}$.  Let $\Lines$ be the moduli space of
  pseudoholomorphic spheres in that same homology class as the given
  sphere (call these \emph{lines}).  By intersection theory arguments
  presented in \cite{McKay:2003} (see below for a little more detail),
  any two lines intersect transversely in a unique point, and (looking
  at the linearized equations) the moduli space $\Lines$ is a smooth
  4-manifold.

  By results of Taubes (see Lalonde \& McDuff
  \cite{Lalonde/McDuff:1996}), there is a unique symplectic structure
  on $\CP{2}$ up to symplectomorphism and rescaling.  Therefore we can
  ensure (possibly by reorienting) that our symplectic structure is
  the usual one. As proven in Gluck \& Warner \cite{Gluck/Warner:1983}
  (and see McKay \cite{McKay:2003} for more details), the space of all
  pseudocomplex structures tamed by a given symplectic structure is
  contractible; simply put this is just the statement that the space
  of strictly contracting maps of a sphere to a hemisphere is
  contractible.  Therefore we can produce a homotopy through smooth
  pseudocomplex structures, which starts at the standard pseudocomplex
  structure, and ends at the given one. Fixing any two points, take
  the line between them in the standard complex structure.

  As proven in Duistermaat \cite{Duistermaat:1972}\footnote{As the reviewer points
out, it is remarkable how much effort the pseudoholomorphic community
could have saved had they been aware of Duistermaat's beautiful work earlier.
At this late stage, Duistermaat's results have probably all been rediscovered
by researchers in pseudoholomorphic curves,
but there is no better place to read them than the original \cite{Duistermaat:1972}.}
(and again see
  McKay \cite{McKay:2003} for more), cokernel computations proceed as
  in the standard structure to show that the problem of constructing a
  pseudoholomorphic sphere through two distinct points is a well-posed
  elliptic system whose linearization has vanishing kernel and
  cokernel, and that as we pass through the homotopy, we can deform
  the sphere to continue to pass through the points.  To give a little
  more detail, $T_{\lambda} \Lines = \Cohom{0}{\nu_{\lambda}}$ where
  $\nu_{\lambda}$ is the normal bundle of $\bar{\lambda}$, equipped
  with the Duistermaat complex structure (see Duistermaat
  \cite{Duistermaat:1972}).  But $\CP{1}$ has only one complex
  structure up to diffeomorphism, and the normal bundle is
  topologically determined to be $\OO{1}$.  To ask that a section of
  this bundle vanish at 2 points forces it to vanish everywhere.
  Moreover, the cokernel is identified with
  $\Cohom{1}{\nu_{\lambda}}=0$, as Duistermaat shows.  Therefore there
  is a pseudoholomorphic sphere of self-intersection 1 through any
  pair of points. Uniqueness is clear by intersection theory, and
  these must be lines. As above, distinct lines intersect at a
  unique point transversely. Therefore the 4-manifold is a projective
  plane.

  Let $\Cor \subset \Pts \times \Lines$ be the incidence
  correspondence. We need to show that it is a smooth embedded
  submanifold of dimension $6$, and that the maps $\Cor \to \Pts$ and
  $\Cor \to \Lines$ are submersions.  By the same ellipticity argument
  that shows that $\Lines$ is a 4-manifold, we find that $\Cor$ is a
  6-manifold, the moduli space of pointed pseudoholomorphic spheres in
  the given homology class.  Suppose that $E \to \Gro{2}{T \Pts}$ is
  our pseudocomplex structure. We know by definition that $E \to \Pts$
  is a fiber bundle with compact fibers. We can map $\Cor \to E$ by
  taking a pointed line to the tangent plane of that line at that
  point. This map is injective and smooth.

  By the same homotopy argument we used above for pairs of points,
  adapted now to tangent planes, we can take any point of $E$, i.e.  a
  point of $\Pts$ and a potential tangent direction $P \in
  \Gro{2}{\Pts}$ for a plane curve, and find a unique line tangent to
  $P$. Therefore $\Cor \to E$ is onto.  The line depends smoothly on
  the choice of tangent plane $P$, by bootstrapping, so $\Cor \to E$
  is a diffeomorphism, and $\Cor \to \Pts$ is a submersion.  We still
  have to show that $\Cor \to \Lines$ is a submersion and that $\Cor
  \to \Pts \times \Lines$ is a embedding. Smoothness of each of these
  maps follows from the previous remarks.

  Thinking of $\Cor$ as a (smooth) moduli space of pointed lines,
  elliptic theory (again looking at the linearized equations) tells us
  that $T_{(p,\lambda)} \Cor$ is the set of pairs
  $\left(\dot{p},\dot{\lambda}\right)$ where $\dot{p} \in T_p \Pts$
  and $\dot{\lambda} \in \Cohom{0}{\nu_{\lambda}}$, so that
  $\dot{\lambda}(p)=\dot{p}$ modulo $T_p \bar{\lambda}$ (there are no
  other obstructions, because those would live in
  $\Cohom{1}{\nu_{\lambda}}=0$). For the same reason, there are no
  base points for this line bundle, since it is just $\OO{1}$, and so
  the tangent space to $\Cor$ has $6$ real dimensions. Indeed it has
  an almost complex structure. Clearly $\Cor \to \Lines$ is a fiber
  bundle.

  Finally, we need to prove that point pencils meet transversely in
  $\Lines$.  This follows immediately from section~6 of
  McKay~\cite{McKay:2003}, where we constructed the dual pseudocomplex
  structure, and from intersection theory presented in the same paper.
  Now we can apply theorem~\vref{thm:BI}.

  In case we start with the manifold $\Pts=\CP{2}$ equipped with a
  pseudocomplex structure tamed by a symplectic structure, as above we
  can assume that the symplectic structure is the Fubini--Study
  symplectic structure, and find a homotopy through pseudocomplex
  structures to the standard complex structure, and homotope lines, to
  ensure that there are pseudoholomorphic spheres in the usual
  homology class.
\end{proof}
\begin{remark}
  Using the results of Gluck \& Warner \cite{Gluck/Warner:1983} (also
  see McKay \cite{McKay:2003} p. 258), we can identify every regular
  4-dimensional projective plane with a smooth fiber bundle map
  $\Sminus \to \Splus$ (between the unit sphere bundles of the bundles
  of anti-self-dual and self-dual 2-forms) which are strictly
  contracting on each fiber, and have image contained in the same
  hemisphere of $\Splus$ that contains the Fubini--Study symplectic
  form. This description gives the differential system explicitly, and
  makes the taming symplectic structure manifest, but leaves the lines
  as unknown solutions of a differential system. It is also not
  functorial, since families might have varying symplectic structures.
  The point of view of a 4-dimensional projective plane as a map $\Cor
  \to \Pts \times \Lines$ makes the lines explicit, but leaves the
  symplectic structure hidden.  Either description parameterizes the
  regular 4-dimensional projective planes with 2 functions of 6
  variables.
\end{remark}
\begin{proposition}
  Every elliptic system of 2 equations for 2 functions of 2 variables
  occurs locally as the differential system for plane curves on a
  regular 4-dimensional projective plane.
\end{proposition}
\begin{proof}
  The system can be locally symplectically tamed (see McKay
  \cite{McKay:2003}).  In local Darboux coordinates, take an open set
  of small volume, and paste it into $\CP{2}$ matching up with the
  usual Fubini--Study symplectic structure on $\CP{2}$. The picture of
  Gluck \& Warner shows us that we can glue together, outside of some
  compact set, this elliptic equation with the usual one for complex
  curves in $\CP{2}$.
\end{proof}
\section{Smooth but irregular 4-dimensional projective planes}
\begin{remark}
  This last approach might provide a mechanism to construct irregular
  smooth projective planes. We might hope to take a regular projective
  plane, thought of in the Gluck \& Warner picture as a strictly
  contracting bundle map $\Sminus \to \Splus$, and deform it while keeping it
  symplectically tamed into a map which is contracting but not
  strictly. One has to prove that the lines survive as smooth
  surfaces, in a smooth family, forming a smooth projective plane. The
  lines would remain symplectically tamed, so this \emph{a priori}
  estimate, together with some local analysis of degenerations of
  elliptic systems, might build examples of irregular smooth
  projective planes.
\end{remark}
\begin{lemma}
  Let $\Pts$ be a smooth 4-dimensional projective plane.  Identify
  $\Pts=\CP{2}$ by a symplectomorphism (whose existence is ensured by
  theorem~\vref{thm:CPTwo}). Take the usual Fubini--Study metric on
  $\CP{2}$.  Note that this makes the symplectic form a self-dual
  2-form. Identify 2-planes in tangent spaces of $T_p \CP{2}$ with
  unit 2-forms in $\Lm{2}{T^*_p \CP{2}}$ by the Pl{\"u}cker embedding:
  \[
  \operatorname{Pl} : \Gro{2}{T \CP{2}} \to \Lm{2}{T^* \CP{2}}.
  \]
  Let $\Splus, \Sminus \subset \Lm{2}{T^* \CP{2}}$ be the bundles of unit
  length self-dual and anti-self-dual 2-forms. Then under the Gauss
  map, followed by the Pl{\"u}cker embedding, the space $\Cor$ is
  smoothly homeomorphically mapped to the graph of a fiber bundle
  morphism $\Sminus \to \Splus$ which is contracting on each fiber. The image
  of each fiber lies inside the hemisphere containing the symplectic
  form. The map is strictly contracting just when the projective plane
  is regular, and is $C^{k-1}$ if the projective plane is $C^k$.
  Continuous isomorphisms of smooth projective planes determine
  continuous isomorphisms of the maps $\Sminus \to \Splus$.
\end{lemma}
\begin{proof}
  For regular 4-dimensional projective planes, see McKay
  \cite{McKay:2003}, section 2.6, where it is shown that the map $\Sminus
  \to \Splus$ is strictly contracting.  Consider an irregular 4-dimensional
  projective plane. By corollary~\vref{cor:generic}, the generic
  point is regular.  So locally the image of $\operatorname{Pl} g$
  is locally the graph of a strictly contracting map near regular
  points. By continuity (see theorem~\vref{thm:BoediTwo}), the map
  thus locally defined is contracting at every point. However, it
  is not clear that this subset of $\Sminus \times \Splus$ is globally
  the graph of a map. By the topological pigeonhole
  principle, it is enough to show that projection of the Pl{\"u}cker map
  to $\Sminus$ is 1-1. We can restrict to studying $\Sminus_p,\Splus_p$
  at a particular point $p \in \Pts$. Lets write $S^{+}$
  for $\Splus_p$ and $S^{-}$ for $\Sminus_p$

  Consider the map $\left(\sigma+,\sigma_-\right)$. By Sard's
  lemma, the set of irregular values has measure 0. By Fubini's
  theorem, the set of points $s_- \in S^{-}$ for which
  the irregular values in $S^+ \times s_-$ have positive measure
  has measure zero.

  In homology, let $A=\left(\sigma_+,\sigma_-\right)\left[\bar{p}\right]
  \in \Homol{2}{S^+ \times S^-}$. The image of $\sigma_+$ is contained
  in the hemisphere around the symplectic form (this precisely expresses
  the condition that the symplectic form tames the smooth projective
  plane; see McKay \cite{McKay:2003} 2.4), so $A$ has no
  component in $\Homol{2}{S^+}$ and so $A=d \left[S^-\right]$,
  some integer $d$. Consider the intersection with the
  graph of an isometry $S^- \to S^+$. For example, the subGrassmannian
  $\grLine{\ell}{2}{T_p \Pts}$ is such a graph (see McKay \cite{McKay:2003}
  2.4). If we pick a real line $\ell \subset T_p \bar{\lambda}$ for a regular
  $\lambda$, then the intersection point is unique, $\lambda=\lambda(\ell)$
  the magnification. But moreover, the intersection is negative
  (see McKay \cite{McKay:2003} 2.4 again). The isometry in
  this case is orientation reversing, and its graph has homology
  $\left[S^+\right]-\left[S^-\right]$ (once again \cite{McKay:2003} 2.4).
  Therefore $d=1$. So the map $\sigma_-$ has degree 1.

  If $\sigma_-$ is not 1-1, then image of the map $\operatorname{Pl} g$
  inside $S^+ \times S^-$ must intersect some submanifold
  $S^+ \times s_-$ at 2 points at least. For every $s_-$, there
  must be an intersection, because $\sigma_-$ is onto. For a full
  measure set of choices of $s_-$, all of the intersection
  points will be regular values of $\left(\sigma_+,\sigma_-\right)$,
  and therefore will be points of positive intersection.
  But the intersection number is 1, so for a full measure open set
  of points $s_- \in S^-$, there will be a unique transverse
  point of intersection. The same will be true if we replace
  $S^+ \times s_-$ with the graph of any strictly contracting
  map $S^+ \to S^-$. This is as far as I can get with
  homology arguments, but it is satisfying to see that we
  can nearly prove the result this way.

  Suppose that $\sigma_-\left(\lambda_1\right) =
  \sigma_-\left(\lambda_2\right)$,
  but $\sigma_+\left(\lambda_1\right) \ne \sigma_+\left(\lambda_2\right)$.
  Write the relevant 2-planes as $\sigma^1_+ + \sigma_-$
  and $\sigma^2_+ + \sigma_-$. Since the Fubini--Study
  form $\omega_+$ tames both 2-planes, we must have
  both $\sigma^j_+$ lying in the same hemisphere as $\omega_+$,
  and therefore $0<\sigma^1_+ \wedge \sigma^2_+<1$.
  Let $1-\epsilon=\sigma^1_+ \wedge \sigma^2_+$.
  Therefore
  \begin{align*}
  \left(\sigma^1_+ + \sigma_-\right)
  \wedge
  \left(\sigma^2_+ + \sigma_-\right)
  &=
  1-\epsilon-1\\
  &=-\epsilon < 0.
  \end{align*}
  But the intersections are positively oriented.
\end{proof}
\begin{corollary}
The differential system for curves on any smooth 4-dimensional
projective plane is a topological submanifold of the Grassmann
bundle, and is the limit in $C^0$ topology of the differential
system for curves in a regular 4-dimensional projective plane.
\end{corollary}
\begin{remark}
The purpose of this corollary is the give evidence for the
conjecture that every $C^k$ smooth projective plane is a limit
of $C^k$ regular projective planes, in the $C^k$ topology.
Note the importance of the Radon transform.
\end{remark}
\begin{remark}
  Take $F_0 \to \Pts \times \Lines_0$ a smooth 4-dimensional
  projective plane, perhaps not regular, and assume $\Pts=\CP{2}$
  bears the standard symplectic structure, positive on the lines of
  $F_0$. Take $\sigma_0 : \Sminus \to \Splus$ the associated Gluck--Warner
  map. From this lemma, we can construct a smooth family $\sigma_t :
  \Sminus \to \Splus$ of strictly contracting maps approaching $\sigma_0$.
  Let $F_t \subset \Gro{2}{T \Pts}$ be the associated embedded
  submanifold (the inverse image under the Pl{\"u}cker map). Take
  $\Lines_t$ the moduli space of lines (i.e. degree 1 maps $\CP{1} \to
  \Pts$, $F_t$-holomorphic).  Let $\widetilde{\Lambda}$ be the set of
  pairs $(t,\lambda)$ where $t \ge 0$ and $\lambda$ is an $F_t$-line.
  Elliptic theory tells us that $\widetilde{\Lambda}$ is a smooth
  cobordism away from the boundary $t=0$. If we could show that
  $\widetilde{\Lambda}$ is a smooth cobordism, then we would see that
  irregular 4-dimensional planes are all limits of regular planes.
\end{remark}
\begin{remark}
  Plane curve theory in irregular smooth projective planes is
  thus a special case of singular perturbation theory of
  first order determined elliptic PDE for two functions of
  two variables.
\end{remark}
\section{Maps taking curves to curves}
\begin{definition}
  A \emph{curve morphism} of smooth projective planes is a map taking
  plane curves (thought of as subsets) to plane
  curves.  For instance, every diffeomorphic collineation is a curve
  morphism.
\end{definition}
\begin{theorem}\label{thm:ContIso}
  A continuous [homeomorphic] curve isomorphism of regular projective
  planes of dimension 8 or 16 [4] is a smooth collineation.  In
  particular, smooth isomorphisms of the differential system for plane
  curves are smooth collineations.
\end{theorem}
\begin{proof}
  In 8 or 16 dimensional regular projective planes, curves are pieces
  of lines, so the isomorphism must take lines to lines, hence a
  collineation. In 4 dimensional regular projective planes similarly,
  any line must be taken to a pseudoholomorphic curve. Homeomorphism
  ensures that it is in the appropriate cohomology class. By the
  intersection theory of McKay \cite{McKay:2003}, it must be a line.
  Therefore the map is a collineation.  By theorem~\vref{thm:BK}, a
  continuous collineation is smooth.
\end{proof}
This strengthens the rather poor results of McKay \cite{McKay:2004}.
It is surprising, since the differential system for plane curves on a
4 dimensional projective plane can be very flexible locally (see McKay
\cite{McKay:2003} for examples, besides the classical $\CP{2}$).
\section{Local characterization of classical projective planes}
\begin{proposition}
  A smooth projective plane of dimension 2 or 4 has all of its tangent
  trialities classical just when it is regular.
\end{proposition}
\begin{proof}
  The proof in dimension 4 is a long calculation; see McKay
  \cite{McKay:2003}.
\end{proof}
\begin{theorem}[Tresse]
  A smooth projective plane of dimension 2 is classical just when the
  Tresse invariants vanish.
\end{theorem}
\begin{proof}
  See Tresse \cite{Tresse:1896}, Arnol{\cprime}d \cite{Arnold:1988},
  Cartan \cite{Cartan:70,Cartan:153}.
\end{proof}
\begin{theorem}
A smooth projective plane of dimension 4 is isomorphic to $\CP{2}$
just when it and its dual are regular, and the differential
invariants $T_2$ and $T_3$ of pseudocomplex structures (presented in
McKay \cite{McKay:2003}) vanish. Equivalently, a smooth projective
plane of dimension 4 is isomorphic to $\CP{2}$ just when it and its
dual bear almost complex structures for which the line and pencils
are pseudoholomorphic curves.
\end{theorem}
\begin{remark}
Note that it is not sufficient to check $T_3$ just for the smooth
projective plane, since the smooth projective planes with $T_3=0$
are precisely the almost complex structures on $\CP{2}$ tamed by the
usual symplectic structure (an open condition), and these form an
infinite dimensional family, depending on 8 functions of 4
variables. Vanishing of $T_2$ alone is not invariantly defined,
since $T_2$ is a relative invariant, and arbitrary multiples of
$T_3$ can be added to it.
\end{remark}
\begin{proof}
As shown in McKay \cite{McKay:2003}, a smooth projective plane is an
almost complex manifold, and its curves pseudoholomorphic curves,
just when $T_3=0$.  In terms of that paper, this is just the
condition that $\tau_1 \wedge \bar{\omega}=0$, which is shown to be
equivalent to $0=\sigma=\tau_1 \wedge \bar{\omega}=\tau_3 \wedge
\bar{\omega}$.  If the dual is also almost complex, then this forces
the same equations on the dual invariants.  McKay \cite{McKay:2003}
section~6 shows that $T_3$ of the dual plane is $T_2$ of the
original plane.  This leaves only the invariants $U_3$ and $V_2$.
Writing out the structure equations of McKay \cite{McKay:2003} pg.
20, and differentiating once, reveals that these are also forced to
vanish.  This forces the differential system for curves to be
isomorphic to the Cauchy--Riemann equations, so by
theorem~\vref{thm:ContIso}, the smooth projective plane is
isomorphic to the classical model $\CP{2}$.
\end{proof}
\begin{remark}
This proof requires $C^5$ differentiability, to define and
differentiate all of the required invariants.
\end{remark}
\begin{remark}
A similar approach via the method of equivalence could certainly
determine local invariants for smooth projective planes of dimension
8 and 16 whose vanishing is necessary and sufficient for isomorphism
with the model projective planes.
\end{remark}
\section{Curves with singularities}
It will be helpful to adopt a more general notion of plane curve in a
regular 4-dimensional projective plane, as in McKay \cite{McKay:2003}
p. 280; the reader will need a copy of that article in hand to follow
our arguments from here on. Essentially the idea is that $\Cor$ admits
a canonical almost complex structure for which $\Theta$ is a complex
subbundle of the tangent bundle, and for which the fibers of $\Cor \to
\Pts$ and $\Cor \to \Lines$ are pseudoholomorphic curves, and for
which generalized plane curves are pseudoholomorphic.  This almost
complex structure is derived in McKay \cite{McKay:2003}. A plane curve
is a pseudoholomorphic map $\Phi : C \to \Cor$ (not necessarily an
immersion) which is tangent to $\Theta$, i.e. for which all 1-forms
$\vartheta$ on $\Cor$ vanishing on $\Theta$ pull back to $\Phi^*
\vartheta=0$.  If such a map is Lipschitz, then it is smooth by
elliptic regularity (see McKay \cite{McKay:2003}).
\begin{theorem}[Micallef \& White]\label{thm:MW}
  Every compact plane curve (perhaps with boundary) in a regular
  4-dimensional projective plane is taken by a Lipschitz homeomorphism
  of a neighborhood of a point to a plane curve in the classical
  4-dimensional projective plane $\CP{2}$.  If the tangents to the
  curve at that point all lie tangent to the same line (for example,
  not a nodal point), then we can make a smooth diffeomorphism instead
  of merely Lipschitz.
\end{theorem}
\begin{remark}
  The crucial idea is that this result holds true at singular points,
  and at intersections, whether transverse or not. Thus the
  intersection theory of plane curves is isomorphic to the
  intersection theory of complex curves. For proof see Micallef \&
  White \cite{Micallef/White:1995} and also see Sikorav
  \cite{Sikorav:2000}.
\end{remark}
\begin{remark}
  The reader might be able to generalize the work of
  Duval \cite{Duval:1999} or of Kharlamov \&
  Kulikov \cite{Kharlamov/Kulikov:2003}
  to plane curves in 4-dimensional smooth projective planes.
\end{remark}
\section{Quadrics}
\begin{definition}
  A \emph{quadric} or \emph{conic} is a closed immersed plane curve $\phi :
  C \to \Pts$, whose degree $\phi_* [C]/\left[\bar{\lambda}\right]$ is
  2, so that no path component of $C$ is mapped to a point.
\end{definition}
\begin{lemma}
  Any quadric in a regular 4-dimensional projective plane is either a
  pair of lines, or is a smooth embedded Riemann surface and contains
  5 points with no three of them on the same line.
\end{lemma}
\begin{proof}
  A quadric is a map $Q \to \Cor$ from a compact Riemann surface $Q$
  which is pseudoholomorphic and tangent to $\Theta$, and which maps
  $Q \to \Cor \to \Pts$ to a degree 2 curve. Split $Q$ into components
  $Q = \coprod Q_{\alpha}$.  Each $Q_{\alpha}$ must map to a point or
  have a positive degree.  No components mapping to points are allowed
  by definition, so each component must have positive degree. Because
  the total degree is 2, either $Q$ has 2 components of degree 1,
  which are therefore lines (by intersection theory with their tangent
  lines), or has one component of degree 2. The singularities of the
  map $Q \to \Cor \to \Pts$ must be diffeomorphic to singularities of
  a plane curve in the classical projective plane $\CP{2}$. The
  selfintersection number at each singularity must be at most 2, by
  the intersection theory of McKay \cite{McKay:2003} proposition~13,
  p.290. Moreover, by the Micallef \& White theorem~\vref{thm:MW}, the
  self-intersections must look like those of an algebraic curve, so to
  have selfintersection number at must 2 at a singularity, it must be
  a double point, i.e. $Q$ is an immersed submanifold. But in that
  case, the local picture of two surfaces intersecting transversely,
  we can take the line tangent to one of those surfaces at the
  intersection point, and it will have intersection number 3 or more.
  This contradicts the degree being 2, so the tangent line must have
  infinite order intersection. By McKay \cite{McKay:2003} theorem~4,
  p. 282, this forces a line to be contained in $Q$. Looking at the
  other tangent line, we find that $Q$ is a union of two lines.

  Therefore either $Q$ is smooth or a union of two lines. Suppose that
  $Q$ is smooth.  If $Q$ contains a line, say $\bar{\lambda}$, then
  take any two points $p_0$ and $p_1$ of $Q$ not on $\bar{\lambda}$,
  and draw the line $\mu=p_0 p_1$. The line $\mu$ strikes $Q$ at
  $p_0,p_1$ and a point of $\bar{\lambda}$. Therefore $\mu$ must also
  have infinite order intersection with $Q$, and again $Q$ is the
  union of $\bar{\mu}$ and $\bar{\lambda}$.

  Therefore either $Q$ is smooth containing no lines, or is a union of
  two lines (possibly a double line, i.e. parameterizing a line twice
  over). If $Q$ is smooth, then no three points of $Q$ lie in a line,
  and therefore any 5 points of $Q$ will do.
\end{proof}
\begin{lemma}
  A smooth quadric in a regular 4-dimensional projective plane is
  diffeomorphic to a 2-sphere.
\end{lemma}
\begin{proof}
  Take a smooth quadric $Q$, a point $q \in Q$. Map $p \in \Pts
  \backslash q \mapsto pq \in \bar{q}$, the Hopf fibration. Restricted
  to $Q \backslash q$, this map has smooth inverse $\lambda \mapsto p$
  where $p$ lies on $\bar{\lambda}$ and on $Q$, defined for all
  $\lambda \in \bar{q}$ with two distinct points of intersection with
  $Q$. But every line has either a pair of distinct points of
  intersection with $Q$, or is tangent to $Q$ (double intersection
  point), since higher intersections are forbidden by homology count.
  We can extend the map to take $q$ to the tangent line at $q$ and
  obtain a bijection $Q \to \bar{q}$.  Away from $q$, this map
  identifies a smooth quadric minus a point with a line minus a point.
  By lemma~\vref{lemma:secant}, if we extend the map to take $q$ to the tangent
  line to $Q$, then it extends to a homeomorphism. Because the 2-sphere
  has a unique smooth structure, $Q$ is diffeomorphic to the 2-sphere.
\end{proof}
\begin{remark}
If the tangent plane to $Q$ at $q$ is regular, then
the map $Q \to \bar{q}$ is an immersion, because
$Q$ is never tangent to second order to any of its tangent
lines. Second order tangency implies intersection number
$3$, as shown in McKay \cite{McKay:2003}. If $Q$ has
irregular tangent line at $q$, then it isn't clear.
\end{remark}
\begin{lemma}
  In the canonical conformal structure (see McKay \cite{McKay:2003}),
  a smooth quadric is biholomorphic to $\CP{1}$, and its normal bundle
  is diffeomorphic to $\OO{4}$.
\end{lemma}
\begin{proof}
  The uniqueness of conformal structure on the sphere is well-known.
  The self-intersection number of a smooth quadric must be 4, by its
  homology, giving the topology of the normal bundle. Apply the
  classification of rank 2 vector bundles on the sphere.
\end{proof}
\begin{lemma}
  The space of smooth quadrics is an oriented smooth real manifold of
  dimension $10$.
\end{lemma}
\begin{proof}
  Consider the normal bundle $\nu$. It comes equipped with the
  linearization of the differential system for plane curves.
  Following Duistermaat \cite{Duistermaat:1972} (or McKay
  \cite{McKay:2003}), we can write the linearized operator
  $Lu=\bar{\partial}u+b\bar{u}$, for sections $u$ of the normal
  bundle, turning the normal bundle into a complex line bundle, and
  the smooth quadric into a complex curve biholomorphic to $\CP{1}$.
  The surjectivity of $L$ is proven in Duistermaat
  \cite{Duistermaat:1972} p.238.  The smoothness of the moduli space
  follows by standard elliptic theory.  We can count the dimension of
  the moduli space as the dimension of the kernel of $L$, by following
  Gromov \& Shubin \cite{Gromov/Shubin:1992}:
\[
\dim \ker L = \ind L + \dim \ker L^t,
\]
 (or just using the standard Riemann--Roch theorem, following
  Duistermaat, as the reviewer points out),
and $L^t$ is the adjoint operator on the dual bundle $\kappa \nu^{-1}$
(where $\kappa$ is the canonical bundle).  But this operator also has
the form $L^t v = \bar{\partial}v+B\bar{v}$ (see Duistermaat
\cite{Duistermaat:1972}).  By the Bers similarity theorem (see Bers
\cite{Bers:1954}), if $v$ lies in the kernel of $L^t$, then $v=e^S V$
for $e^S$ a nowhere zero holomorphic section of a trivial line bundle,
and $V$ a holomorphic section of another line bundle, say
$\mathcal{L}$. Counting Chern classes, $c_1\left(\kappa
  \nu^{-1}\right)=-6$. So $V$ must be a holomorphic section of
$\OO{-6}$, and therefore vanishes, so $v$ does as well. By the
Atiyah--Singer index theorem (or again just the Riemann--Roch theorem)
\begin{align*}
  \dim \ker L &= \ind L \\
  &= 2 c_1(\nu) + 2 \\
  &= 10
\end{align*}
as a real vector space.
\end{proof}
\begin{lemma}
  Pick 5 points $p_1,\dots,p_5$ in a regular 4-dimensional projective
  plane, no 3 of them lying in a line. There is a smooth 1-parameter
  family of regular projective plane structures $\Cor \times \R{} \to
  \Pts \times \Lines$ (write $\Cor_t$ for $\Cor \times
  \left\{t\right\}$) and a smooth family of points $P_j : \R{} \Pts$,
  $j=1,\dots,5$, so that $\Cor_0$ is isomorphic to the classical
  projective plane $\CP{2}$, $\Cor_1$ is the regular 4-dimensional
  plane of our hypothesis, all $\Cor_t$ are tamed by the standard
  symplectic structure, and for each time $t$, no 3 of the points
  $P_1(t),\dots,P_5(t)$ are colinear with respect to $\Cor_t$.
\end{lemma}
\begin{proof}
  Take any smooth deformation $\Cor_t$ from the classical structure to
  the given structure, tamed by the usual symplectic structure,
  guaranteed to exist by theorem~\vref{thm:taubes}.  Let $M$ be the
  set of sextuples $\left(q_1,\dots,q_5,t\right) \in \prod^5 \Pts
  \times \R{}$ for which no three of the $q_j$ are colinear with
  respect to $\Cor_t$. Clearly $M$ is an open subset of $\prod^5 \Pts
  \times \R{}$.  Indeed $M$ is just the leftover part when we remove
  the diagonal loci $q_i=q_j$ and the loci where $q_k \in
  \overline{q_i q_j}$.  The diagonal loci are obviously submanifolds
  of codimension 4. The locus where $q_k \in \overline{q_i q_j}$ can
  be written as $\left(q_i,q_j q_k\right) \in \Cor$.  So this locus is
  the inverse image under the map
\[
\left(q_i,q_j,q_k\right) \mapsto \left(q_i,q_jq_k\right) \in \Pts
\times \Lines
\]
of $\Cor$.
\begin{lemma}
  The map $\left(p_1,p_2\right) \in \Pts \times \Pts \backslash \Delta
  \Pts \mapsto p_1 p_2 \in \Lambda$ has full rank at all points where
  it is defined (i.e.  where $p_1 \ne p_2$).
\end{lemma}
\begin{proof}
  Pick two points $p_1 \ne p_2$, and take two lines
  $\lambda_1,\lambda_2 \ne p_1 p_2$.  Consider the smooth map $\lambda
  \mapsto \left(\lambda \lambda_1, \lambda \lambda_2\right)$ defined
  where $\lambda \ne \lambda_1$ and $\lambda \ne \lambda_2$. This map
  provides a local section $\text{open } \subset \Lambda \to \Pts
  \times \Pts \backslash \Delta \Pts$. Therefore the map has full
  rank.
\end{proof}
Therefore, returning to our map $\left(q_i,q_j,q_k\right) \mapsto
\left(q_i,q_jq_k\right)$, this map has full rank, and therefore the
inverse image of $\Cor$ is a smooth submanifold of codimension $2$.
Our manifold $M$ is the complement in $\prod^5 \Pts \times \R{}$ of a
finite set of codimension 2 and codimension 4 submanifolds, and is
therefore connected.
\end{proof}
\begin{lemma}
  Given any 5 points in a regular 4-dimensional projective plane, no 3
  of which lie on a line, there is a unique smooth quadric through
  them, and the quadric depends smoothly on the choice of the 5
  points.
\end{lemma}
\begin{remark}
  The proof is Gromov's \cite{Gromov:1985} $2.4.B_1''$, with some more
  details.
\end{remark}
\begin{proof}
  Start by deforming the 5 points $P_1(t),\dots,P_5(t)$ and the
  projective plane structure $\Cor_t$ from the classical projective
  plane at $\Cor_0$, to the given projective plane at $\Cor_1$. Ensure
  that the no three of the points $P_j(t)$ lie in a line, at any time
  $t$, and that all of the $\Cor_t$ projective plane structures are
  tamed by the same symplectic structure. Gromov compactness ensures
  that the set of smooth quadrics through the given points is compact
  on any compact interval of $t$ values in $\R{}$. There is a unique
  such quadric at $t=0$, by classical plane algebraic geometry. To
  ensure the survival of the smooth quadric on any open subset of
  $\R{}$, we employ the continuity method of elliptic partial
  differential equations. This reduces to showing the surjectivity of
  the linear operator
\[
Lu=\bar{\partial}u+b\bar{u}
\]
where $\bar{\partial}+b$ is derived in Duistermaat
\cite{Duistermaat:1972}, p. 237, $u$ a section of the normal bundle
$\nu$, and we need surjectivity among sections $u$ with specified
values $u\left(P_j\left(t\right)\right)$ at the 5 points. By
linearity, we can assume that those specified values are $u=0$.
(As the reviewer points out, surjectivity here is
proven in Ivashkovich \& Shevchisin \cite{Ivashkovich/Shevchisin:1999}
and Barraud \cite{Barraud:2000}.)

The surjectivity of this operator is equivalent to the injectivity of
the adjoint operator, which is
\[
v \mapsto \left( \left(\bar{\partial}+b\right)^t,v\left(P_1(t)\right),
  \dots,v\left(P_5(t)\right) \right),
\]
$v$ a section of $\kappa \nu^{-1}$ and $\kappa$ the canonical bundle.
Topologically, the first Chern classes are
\begin{align*}
  c_1\left(\kappa \nu^{-1}\right)&=
  c_1\left(\kappa\right)-c_1\left(\nu\right)\\
  &=-2-4\\
  &=-6.
\end{align*}

Consider the divisor $\mu=P_1(t)+\dots+P_5(t)$.  By the
Gromov--Shubin--Riemann--Roch theorem (see Gromov \& Shubin
\cite{Gromov/Shubin:1992,Gromov/Shubin:1993,Gromov/Shubin:1994}), the
space of solutions of $L$ with zeros on $\mu$ satisfies
\[
\dim \ker (L,\mu) = \ind L - \deg \mu + \dim \ker
\left(L^t,-\mu\right).
\]
The index of $L$, by the Atiyah--Singer index theorem, is $c_1(L)+1$
(see Duistermaat \cite{Duistermaat:1972}, p.238).  The number $\deg
\mu $ depends only on $\mu$, not on the operators involved (see Gromov
\& Shubin \cite{Gromov/Shubin:1994}, p.169), so we can calculate it
for a smooth quadric in the standard projective plane, and find that
$\deg \mu=5$. Therefore
\[
\dim \ker (L,\mu) = \dim \ker \left(L^t,-\mu\right).
\]
Applying the Bers similarity principle (see Bers \cite{Bers:1954}),
any $u$ in $\ker (L,\mu)$ has the form $u=e^sU$ where $e^s$ is a
nowhere vanishing section of a (obviously trivial) line bundle, and
$U$ is a section of a line bundle $\nu'$ with the same topology as
$\nu$. By the Birkhoff--Grothendieck theorem, $\nu'=\OO{4}$.  The
section $U$ will have the same zeros as $u$, and so will have 5 zeros
(at the $P_j(t)$). This forces $U=0$.  Therefore $\ker (L,\mu)=0$, and
so $\ker \left(L^t,-\mu\right)=0$, ensuring that $L$ is surjective.
\end{proof}
\begin{remark}
  We didn't really need to use the full force of Bers's similarity
  principle for line bundles; it is enough to notice that $u$ has only
  positive intersections with the zero section, which follows from the
  local Bers similarity principle.
\end{remark}
\begin{corollary}
  The space of smooth quadrics is connected.
\end{corollary}
\begin{proof}
  The space of 5-tuples with no 3 points colinear is the complement of
  a codimension 2 subset of the space of 5-tuples, and therefore is
  connected, and maps smoothly onto the space of smooth quadrics.
\end{proof}
\begin{lemma}
  Take $Q$ a smooth quadric in a regular 4-dimensional projective
  plane.  The map $q \in Q \to T_q Q \in \Lines$ takes $Q$ to a smooth
  quadric in $\Lines$.
\end{lemma}
\begin{proof}
  By corollary~\vref{cor:dualCurves}, this is the map to the dual
  curve, so it is a plane curve (perhaps with singularities).  By
  intersection theory of McKay \cite{McKay:2003}, no line can be
  tangent to a quadric at two points, so the map is injective. The map
  has full rank, because no smooth quadric can be tangent to higher
  than first order with a line. By deforming to the classical case, we
  can compute the degree. Therefore the dual curve is a quadric. If
  the dual curve is not smooth, then it must be a pair of lines. But
  the dual of a line is a point, so the original curve would have been
  a pair of points. Therefore the dual curve is a smooth quadric.
\end{proof}
\begin{lemma}
  There are precisely 4 lines simultaneously tangent to any pair of
  distinct smooth quadrics, counting multiplicities.  We don't have to
  count multiplicities unless there are one or two points of
  nontransverse intersection.
\end{lemma}
\begin{proof}
  Lines tangent to $Q$ are points of the dual curve $Q^*$.  Count
  intersections of the dual curves, which are smooth quadrics.
\end{proof}
\begin{lemma}
  Given a smooth quadric, the smooth quadrics nowhere tangent to it
  form a dense open subset of the quadrics.
\end{lemma}
\begin{proof}
  Pick 4 distinct points on the given quadric, and one point not on
  the given quadric, not on a line through any two of the 4 points.
  The quadric through those 5 points is nowhere tangent to the
  original quadric, because it has 4 distinct points of intersection,
  so by homology counting and positivity of intersection (see McKay
  \cite{McKay:2003}), the two quadrics are nowhere tangent. Given any
  quadric which is tangent, we can pick our 4 intersection points
  close to its intersection points, and our fifth point close to it.
\end{proof}
\begin{lemma}
  Given any pair of nowhere tangent smooth quadrics, we can deform the
  projective plane structure $\Cor \to \Pts \times \Lines$ into the
  classical one, and deform the quadrics so that they remain quadrics
  throughout the deformation, and keep them from every becoming
  tangent.
\end{lemma}
\begin{proof}
  Given one smooth quadric $Q_1$, pick any 4 distinct points on it,
  $p_1,\dots,p_4 \in Q_1$. Draw the tangent line $\lambda$ to $Q_1$ at
  $p_1$, and pick any point $p_5$ of $\bar{\lambda}$ other than
\[
p_1 ,\lambda\left(p_ip_i\right), i,j=1,\dots,4.
\]
Then the quadric through the points $p_1,\dots,p_5$ is smooth and
nowhere tangent to $Q_1$. The space of choices of the $p_j$ points is
clearly a connected manifold, a 4!-fold covering space of the space of
ordered pairs of nowhere tangent quadrics, of dimension $20$, and a
fiber bundle over the space of smooth quadrics. We can easily add a
parameter $t$ to the construction, and look at points
$p_1(t),\dots,p_5(t)$, and not lose the connectedness.
\end{proof}
\begin{lemma}\label{lemma:involution}
  Take $Q$ a smooth quadric, $X$ the set of pairs $(p,\lambda)$ for
  which $p \in Q$ and $\lambda$ is a line through $p$. Map $\iota :
  (p,\lambda)\in X \to \left(p',\lambda\right) \in X$ where $p,p'$ are
  the points where $\bar{\lambda}$ intersects $Q$, and take $p' \ne p$
  if possible, i.e. unless $\bar{\lambda}$ is tangent to $Q$ at $p$.
  The map $\iota : X \to X$ is a smooth diffeomorphism.
\end{lemma}
\begin{proof}
  Clearly $X = \pi_{\Pts}^{-1} Q$ is a smooth manifold.  Where
  $\lambda$ is not tangent to $X$, we can ensure the result by
  transversality. For tangent $\lambda$, the result is immediate in
  Micallef--White coordinates (see theorem~\vref{thm:MW}).
\end{proof}
\begin{remark}
  Pascal's mystic hexagon apparently does \emph{not} give a mechanism for drawing smooth
  quadrics; see Hofmann \cite{Hofmann:1971}.
\end{remark}
\begin{remark}
  Gromov \cite{Gromov:1985} p. 309 $0.2.B$ suggests that the approach
  we have taken here to construct quadrics can construct plane curves
  of all genera. The details of the argument have never been provided;
  Gromov (p. 338 $2.4.B_1''$) suggests that there are some subtleties.
\end{remark}
\section{Poncelet's porism}
Because regular 4-dimensional projective planes are symplectomorphic
to $\CP{2}$, they share the same Gromov--Witten invariants, so that a
huge collection of enumerative problems about plane curves have the
same solutions.  Lets consider some plane geometry which is not
enumerative.  We will search for an analogue of the elliptic curve in
the classical proof of \emph{Poncelet's porism}.  For proof in the
classical 4-dimensional projective plane, see Griffiths \& Harris
\cite{Griffiths/Harris:1977,Griffiths/Harris:1978} and Schwartz
\cite{Schwartz:2001}.
\begin{definition}
  A \emph{polygon} in a smooth projective plane is an ordered
  collection of distinct points $p_1,\dots,p_n,p_{n+1}=p_1$.  The
  lines $p_j p_{j+1}$ are called the \emph{edges} of the polygon,
  while the points $p_j$ are called the \emph{vertices}. A polygon is
  \emph{circumscribed} about a quadric if every edge of the polygon is
  tangent to the quadric. A polygon is \emph{inscribed} in a quadric
  if every vertex lies in the quadric. Given two quadrics, a polygon
  circumscribed about the first one, and inscribed in the second, is
  called a \emph{Poncelet polygon} of those quadrics.
\end{definition}
The classical theorem:
\begin{theorem}[Poncelet]
  For a given pair of nowhere tangent quadrics in the classical
  projective plane of dimension 4 or more, every Poncelet polygon of
  those quadrics belongs to a smooth family of distinct Poncelet
  polygons of those same quadrics, with any one of the vertices [or
  edges] being drawn over the entire quadric in which the polygon
  remains inscribed [circumscribed].
\end{theorem}
\begin{remark}
  For smooth projective planes of dimension 8 or 16, the result is a
  triviality, since quadrics are pairs of lines. Henceforth, consider
  a regular projective plane of dimension 4. It seems very unlikely
  that the Poncelet porism is true for generic regular 4-dimensional
  projective planes.
\end{remark}
Consider a Poncelet polygon. Let $Q_E$ be the quadric that the polygon
circumscribes, and $Q_V$ the quadric in which the polygon is
inscribed.
\begin{lemma}
  Let $T$ be the set of pairs $(p,\lambda)$ so that $p \in Q_V$ and
  $\lambda$ is a line containing $p$ and tangent to $Q_E$.  Map $T \to
  Q_E$, by taking each line to its point of intersection with $Q_E$.
  This map is well-defined, smooth, and a double covering branched at
  $4$ points.
\end{lemma}
\begin{proof}
  Let $\delta : Q_E \to Q_E^*$ be the duality map, taking each point
  to its tangent line.  The map $(1,\delta):Q_V \times Q_E \to Q_V
  \times Q_E^*$ identifies $T$ with the set $T'$ of pairs $(p,q) \in
  Q_E \times Q_V$ so that either $p=q$ and $Q_E$ is tangent to $Q_V$
  at $p$ or $p \ne q$ and $pq=\delta(q)$. Since the map $(p,q)\mapsto
  pq$ has full rank, $T'$ is a submanifold of codimension $2$, except
  possibly at points of the form $(p,p)$, i.e. tangent points of the
  two smooth quadrics.

  Deform to the classical case, deforming the two smooth quadrics
  $Q_E$ and $Q_V$, and forgetting about the Poncelet polygon for the
  moment, and you see that $T$ deforms into the elliptic curve of the
  classical case, so $T$ is diffeomorphic to a 2-torus.  The map $T
  \to Q_E$ is just the composition $T \to Q_E^* \to Q_E$ of the
  projection with the dual map, therefore a smooth map.
\end{proof}
We have two involutions defined on $Q_V \times Q_E^*$:
$\iota_{\Pts}(p,\lambda)=\left(p',\lambda\right)$ where $p,p'$ are the
points of intersection of $\lambda$ with $Q_V$, and
$\iota_{\Lines}(p,\lambda)=\left(p,\lambda'\right)$, where
$\lambda,\lambda'$ are the points of $Q_E^*$ which contain $p$. These
maps are well-defined, except where the intersections are double
points, where we take $p=p'$ ($\lambda=\lambda'$ respectively).
\begin{lemma}
  The maps $\iota_{\Pts},\iota_{\Lines}$ are smooth, and each has 4
  fixed points.
\end{lemma}
\begin{proof}
  The fixed points of $\iota_{\Pts}$ are obviously the points
  $(p,\lambda)$ where $\lambda$ is tangent to both $Q_V$ and $Q_E$
  (i.e. in $Q_V^* \cap Q_E^*$). Dually, the fixed points of
  $\iota_{\Lines}$ are the points $(p,\lambda)$ where $p$ belongs to
  both $Q_V$ and $Q_E^*$.  By transversality, both maps are smooth
  away from their fixed points.  Lemma~\vref{lemma:involution} assures
  smoothness near tangent points.
\end{proof}
\begin{lemma}
  The map $\varpi = \iota_{\Lambda} \iota_{\Pts}$ is a diffeomorphism
  isotopic to the identity.  Neither $\varpi$ nor $\varpi \circ
  \varpi$ have any fixed points.
\end{lemma}
\begin{proof}
  The fixed points of $\varpi$ are the points $(p,\lambda)$ for which
  the line $\lambda$ has double intersection with $Q_E$ at $p$, so
  $\lambda$ is tangent to $Q_E$ at $p$, and for which there is only
  one line, $\lambda$, tangent to $Q_V$ passing through $p$, so
  $\bar{p}$ has double intersection with $Q_V^*$ at $\lambda$, and
  therefore so $\bar{p}$ is tangent to $Q_V^*$ at $\lambda$. But the
  tangent line to $Q_V^*$ at $\lambda$ is $\bar{q}$ for $q \in Q_V$
  the corresponding point, so $\bar{p}=\bar{q}$, i.e. $p=q$.
  Therefore $p \in Q_E \cap Q_V$.  So $p$ is one of the four points of
  $Q_V \cap Q_E$, and $\lambda$ is tangent to both $Q_V$ and $Q_V$,
  necessarily at $p$. Therefore $Q_V$ and $Q_E$ have a common tangent,
  contradicting our hypotheses.

  Consider a fixed point of $\varpi^2$, i.e.
  $\varpi(p,\lambda)=\left(p',\lambda'\right)$ and
  $\varpi\left(p',\lambda'\right)=\left(p,\lambda\right).$ Then
  $\lambda$ and $\lambda'$ are two tangent lines to $Q_E$ which both
  intersection $Q_V$ at both points $p$ and $p'$. Therefore $p=p'$ or
  $\lambda=\lambda'$. If $p=p'$, then $\lambda$ and $\lambda'$ are
  both tangent to both $Q_V$ and $Q_E$, and both must be tangent to
  $Q_V$ at $p=p'$, so must be equal. Dually, if $\lambda=\lambda'$
  then $p=p'$.  Therefore $\varpi(p,\lambda)=(p,\lambda)$,
  contradicting the last paragraph.  Under deformation to the
  classical projective plane, $\varpi$ is taken by isotopy to a
  translation on an elliptic curve, and therefore is isotopic to the
  identity map.
\end{proof}
\begin{lemma}
  Poncelet polygons are precisely periodic orbits of $\varpi$.
\end{lemma}
\begin{proof}
  $\varpi$ by definition takes a point of $Q_V$ and tangent line to
  $Q_E$ to another such point and line, with the next point also
  contained in the first line.  Thus a periodic orbit of $\varpi$
  draws a closed Poncelet polygon.
\end{proof}
\begin{remark}
  Given any smooth quadric, and 5 tangent lines to it
  $\lambda_1,\dots,\lambda_5$, their pairwise intersections
  $\lambda_i\lambda_{i+1}$ lie on a quadric. For a generic choice of 5
  tangent lines, one should be able to show that the resulting quadric
  is smooth, so that there should be Poncelet pentagons. In a generic
  smooth projective plane, these might be the only Poncelet polygons.
\end{remark}
\begin{remark}
  Take an orbit of $\varpi$, and pick an oriented real surface
  $\Sigma$, say a smooth $CW$ complex.  Then average its number of
  intersections with the first $k$ tangent lines from that orbit. Let
  $k$ get large.  I imagine that either a symplectic form emerges, in
  the cohomology class of the usual Radon--Fubini--Study type
  symplectic forms, or a Poncelet polygon, but I have no proof.
\end{remark}
\bibliographystyle{amsplain} \bibliography{poncelet}
\end{document}